\documentclass{gtpart}

\usepackage{pinlabel}
\usepackage{color}

\title{Pseudoholomorphic quilts and Khovanov homology}

\author{Reza Rezazadegan} 
\givenname{Reza}
\surname{Rezazadegan}
\address{Laboratoire de Math\'{e}matiques Jean Leray\\ Universit\'e de Nantes  
         \\France}
\email{reza.rezazadegan@univ-nantes.fr}

\arxivreference{}  
\arxivpassword{}   

\usepackage{khodesh}

\begin{document}


\begin{abstract}
We further study the symplectic Khovanov homology of Seidel and Smith
and its generalization to even tangles.
This homology theory is a conjectural geometric model for Khovanov
homology. In this paper we uncover structures on symplectic Khovanov
homology which have analogues in Khovanov homology.

To each elementary (as well as minimal) cobordism between two tangles
we associate a homomorphism between the symplectic Khovanov homology
groups of the two tangles.
We define the symplectic analogues $H_{s}^{m}$ of Khovanov's arc
algebras  and equip the symplectic Khovanov homology of an $(m,n)$-tangle with
the structure of an $(H_{s}^m,H_{s}^n)$-bimodule. We show that
$H_{s}^{m}$ and Khovanov's $H^{m}$ are isomorphic as associative
algebras over $\Z/2$.  We also obtain a skein exact triangle for symplectic
Khovanov homology which resembles the one for Khovanov homology.
\end{abstract}

\maketitle


\section{Introduction}


This paper is an addendum to \cite{RR1} and its purpose is to prove some natural properties of the tangle invariant discussed therein. 
The tangle invariant is a generalization of the link invariant developed in \cite{SS} by Seidel and Smith as a candidate for a geometric model for Khovanov homology. It is called \emph{symplectic Khovanov homology} and we denote it by $\khs$. 
Roughly speaking the $\khs$ of a link $K$, given as the closure of a braid $\beta\in Br_{m}$, is defined as the Lagrangian Floer homology of a canonically defined Lagrangian $L$ and its image under a symplectomorphism associated to $\beta$. The Lagrangian $L$  lives in a symplectic manifold diffeomorphic to an $(m,m)$-type Springer fiber.  One can hope that the geometric nature of $\khs$ can provide more insight which is not available for Khovanov homology itself. (See \cite{SS2} for a result in this direction.)

The invariant $\khs$ was generalized to even tangles in
\cite{RR1}. One motivation for this generalization was to make the
proof of the equivalence of $\khs$ and Khovanov homology more feasible
by decomposing a given link into elementary tangles. As an example we
obtained a K\"unneth-type formula for symplectic Khovanov homology of
the unlinked union of two links \cite[section 4.4]{RR1}.

In Section \ref{minimalcobo}  we obtain (elementary) cobordism maps for  $\khs$ i.e. if two tangles $T, T'$ are related by an elementary cobordism $S$, we obtain a homomorphism
\bq\label{intro-cobo}
\khs(S): \khs(T)\ra \khs(T').
\eq
These maps are given by counting \pse quilts.
 We do not attempt at showing that our cobordism maps give a well-defined map when a general cobordism is decomposed into elementary ones in different ways.

In the rest of Section \ref{ResultsChapter} we use the cobordism maps \eqref{intro-cobo} to show that $\khs$ shares some of the properties of Khovanov's invariant.
First of all we define a family of associative algebras $\{H^{m}_{s}\}_{m\in \N}$
which are the symplectic analogues of Khovanov's arc algebras $H^{m}$ from  \cite{functorvalued}. 
We set
\bq 
H^m_{s}:=\khs(id_m)
\eq
 with  product given by the maps associated to minimal cobordisms  from Section \ref{minimalcobo} which in turn are given by counting pseudoholomorphic quilted triangles.
 We show that the rings $H^{m}$ and $H^{m}_{s}$ are isomorphic  as
 algebras over $\Z/2$ (Theorem \ref{HsympH}). To do this we compute
 the Khovanov homology of elementary tangles (Prop. \ref{khofbraid}).
 We equip the abelian group $\khs(T)$, assigned to an $(m,n)$-tangle $T$,  with the structure of  an $(H_{s}^m,H_{s}^n)$-bimodule.

 We provide further evidence for the equivalence of this invariant with Khovanov's combinatorially defined invariant by showing the equivalence for flat (crossingless) tangles and elementary cobordisms between them 
 (Propositions \ref{khs=khflat} and \ref{coboequflat} ). 
In \ref{vanishofdiff}  we prove that for flat tangles Floer data can be chosen in such a way that Floer differential vanishes.

At last in Section \ref{exacttri} we use the exact triangle for
fibred Dehn twists from \cite{WWtriangle} to prove an skein exact triangle for $\khs$ (Corollary \ref{thecone}) which resembles that of Khovanov homology after the collapse of bigrading. An argument similar to the one used by Manolescu and Ozsvath \cite{qalternating}  can  then be used to show that the two invariants (i.e. $\khs$ and bigrading-collapsed Khovanov homology) agree on quasi-alternating links. However we do not  present a  proof of this last claim.

We point out that in almost the same time that  the first draft of this paper appeared online, similar cobordism maps were introduced in an independent work by Jack Waldron \cite{Waldron1}.
 Waldron's cobordism maps are  defined in the original setting of the link invariant of \cite{SS}
 and are given using    
the relative invariants associated to surfaces with strip-like ends by Seidel \cite{seideltriangle}. 
 But our maps are defined using the formalism of quilts \cite{WW} and
 are defined after decomposing the involved tangles  into elementary
 ones. Waldron also shows that his maps are independent (up to sign)
 of the decomposition of the cobordisms.
 We expect that his cobordism maps, when restricted to elementary cobordisms, agree with our maps, when the latter is restricted to cobordisms of links.

\begin{remark}
Symplectic Khovanov homology as defined by Seidel and Smith is an invariant of links given as closure of braids therefore comparison of this invariant and Khovanov homology is nontrivial. See ``symplectic Khovanov homology of crossing diagrams'' in Section 5.6 in \cite{Waldron1}. However the the extension from \cite{RR1} which we use here is assigned to tangle decompositions just as Khovanov's invariant. 
\end{remark}

In the  sections 2 to 4 we review, respectively, Khovanov's invariant of tangles, Seidel-Smith invariant for links and the generalization of the latter to tangles.
The notation we use here for the symplectic invariant is a bit different from that of \cite{RR1}. In that paper we used the notation $\mathpzc{Kh}_{symp}$ but here we use $\khs$ (and $\cf$ for the chain complex).

For the results of this paper to hold for Floer chain complexes with coefficients in $\Z$, one needs coherent orientations on the moduli spaces of pseudoholomorphic quilts used. This relies on work in progress by Wehrheim and Woodward \cite{WWorient}.\\

\textbf{Acknowledgements.} I would like to thank Ivan Smith for pointing out a problem with signs in my original proof of Theorem \ref{HsympH}.

\section{Khovanov homology of even tangles}\label{KhChapter}

Here we recall the basic facts about the Khovanov homology of even tangles.

\subsection{Tangles}

This subsection is a repetition of Section 4.1 in \cite{RR1}.
A tangle $T$ is defined to be a compact one-dimensional submanifold of (a diffeomorphic image of) $\mathbb{C}\times [0,1]$ such that $i(T):=T\cap (\mathbb{C}\times \{0\})\subset \mathbb{R} \times  \{0\}$ and $t(T):= T\cap (\mathbb{C}\times \{1\})\subset \mathbb{R} \times  \{1\}$ and both sets are finite. The second assumption makes $i(T)$ and $t(T)$ ordered sets.  In this thesis we deal only with tangles with an even number of initial points and end points. Such tangles are called \textit{even tangles}.
If $\# i(T)=2m, \#t(T)=2n$ we say $T$ is an $(m,n)$-tangle and write  $m T n$. We also allow $m$ and/or $n$ to be zero.  
\begin{definition}
Two tangles $T, T'$ are called \textit{equivalent} if there is a continuous family ${T_t}$  of tangles for ${t\in [0,1]}$ such that $T_0=T$ and $T_1=T'$ and the order of $i(T_t)$ and of $t(T_t)$ is fixed.
\end{definition}

 \begin{figure}
\centering
\includegraphics[scale=0.6]{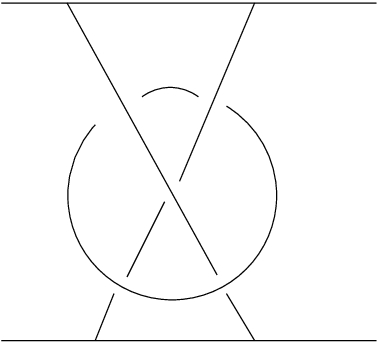}
\caption{A $(1,1)$-tangle}
\end{figure}

\begin{figure}[ht]
\centering
\includegraphics[scale=0.4]{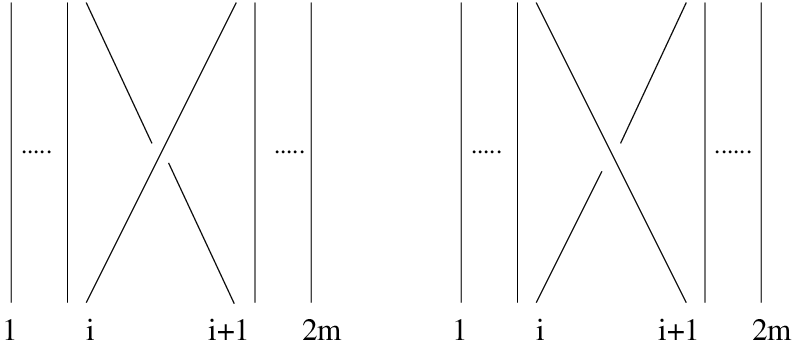}
\caption{The braids $\sigma_i$ and ${\sigma_i}^{t}$}
\label{elembraids}
\end{figure}

\begin{figure}[ht]
\centering
\includegraphics[width=4.5in]{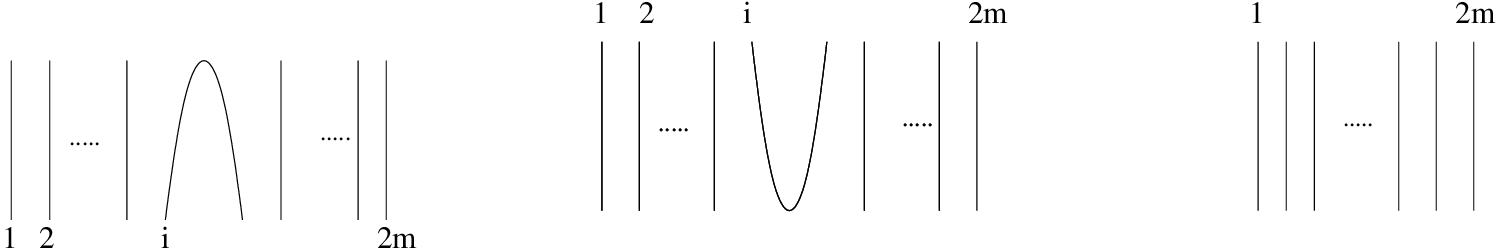}
\caption{$\cap_{i;m}$, $\cup_{i;m}$ and $id_m$}
\label{elementangles}
\end{figure}

Two tangles  $T_1, T_2$ can be composed (concatenated) if $t(T_1)=i(T_2)$. Two equivalence classes $[T_1]$ and $[T_2]$ of tangles can be composed if
$\#t(T_1)=\#i(T_2)$ and composition is done using the ordering on $t(T_1)$ and $i(T_2)$. Composition of tangles is denoted by juxtaposition.
We will use the  notation $id_m$, $\cap_{i;m}$, $\cup_{i;m}$ and $\sigma_{i;m}$ for the elementary tangles in Figures \ref{elembraids} and \ref{elementangles} where $m$ denotes the number of the strands. We might ignore $m$ when there is no confusion. 
When we say a tangle $T$ is equivalent to, say, $\cap_{i;m}$, we implicitly have a one to one correspondence between $i(T)$ and $\{1, 2,\ldots, 2m-2  \}$ and also between $t(T)$ and $\{1, 2,\ldots, 2m\}$ in mind.\\

A \textit{decomposition} of $T$ is a sequence of tangles
\bq\label{tangdeco}n_0 T_1 n_1 T_2\ldots n_{l-1} T_l n_l \qquad n_0=m, n_l=n\eq such that $T$ is equivalent to $T_1 T_2 \cdots T_l$.  
A Morse-theoretic argument shows that  any $T$ can be expressed (not uniquely) as a composition of elementary tangles.
\textit{Crossingless matchings} (section \ref{crossless}) are a special class of $(0,n)$ or $(n,0)$-tangles.  Given a set of $2n$ points on the plane, a crossingless matching is a set of $n$ non-intersecting curves each joining a pair of the given points. In the context of tangles a crossingless matching is regarded as a subset of $\bb{C}\times [0,1]$.
 \begin{definition}
 Let $\pz{C}_n $ be the set of isotopy (in $\mathbb{C}$) classes  of crossingless matchings between $2n$ points on the real line all of whose arcs lie in the upper half plane.
 \end{definition}
 The cardinality of $\pz{C}_n$ equals the $n$'th  Catalan number.

  One can define a category $\mathbf{Tang}$ whose objects are natural numbers and $\hom(m,n)$ consists of equivalence classes of $(m,n)$-tangles. 
  $\mathbf{Tang}$ has a monoidal structure given by putting two tangles $kTl$ and $mTn$ ``side-by-side'' to  obtain a $(k+m,l+n)$-tangle. We denote this by $T\oplus T'$. To each $(m,n)$-tangle $T$ there is assigned a ``mirror image'' $T^t$ which is a $(n,m)$-tangle.
   There is a generators and relations description of $\mathbf{Tang}$ due to Yetter \cite{Yetter} whose proof relies on Reidemeister's description of plane diagram moves.


\begin{lemma}[Yetter \cite{Yetter}]\label{yetterlemma}
The following  are all the commutation relations between elementary tangles where ``$=$'' means equivalence. If $|i-j|>1$ we have:\\\begin{eqnarray}
\sigma_i\sigma_j&=&\sigma_j\sigma_i\\
\cap_{i}\cup_{j}&=&\cup_{j}\cap_{i}\\
\cap_{i}\sigma_{j}=\sigma_{j}\cap_{i}  &\quad& \cup_{i}\sigma_{j}=\sigma_{j}\cup_{i},
\end{eqnarray}
and for any $i$ we have:\\
  \begin{eqnarray}
\sigma_{i}\cup_{i}&=&\cup_{i}\\
\sigma_i \sigma_i^t&=&id\\
\sigma_i\sigma_{i+1}\sigma_i&=&\sigma_{i+1}\sigma_i\sigma_{i+1}\\
\label{capcup}\cap_{i;m}\cup_{i+1;m}&=&id_{m-1} \\
\sigma_i \cup_{i+1} =\sigma_{i+1}^t\cup_i &\qquad& \sigma_i^t \cup_{i+1}=\sigma_{i+1}\cup_i.
\end{eqnarray}
\end{lemma}

\begin{figure}[ht]
\centering
\includegraphics[scale=0.6]{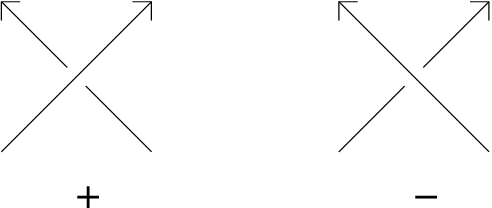}
\caption{Oriented braids $\sig^+$ and $\sig^-$}
\label{posnegbraid}
\end{figure}

The invariants that we consider here are  invariants of \textit{oriented tangles}. An oriented tangle comes with an orientation of each one of its components.  Two example are shown in Figure \ref{posnegbraid}. 

\subsection{The TQFT}\label{khovtqft}
Khovanov homology is based on a 1+1 dimensional TQFT $\mathcal{F}$ whose definition we review here.
1+1 dimensional TQFTs are in one-to-one correspondence with Frobenius algebras. 
 Khovanov \cite{categorification} defines the Frobenius algebra $\pz{V}$ to be equal to $H^*(S^2)\{-1\}$  (i.e. the cohomology of $S^2$ with its grading shifted down by one) as a ring. Let $\pz{1},\pz{X}$ be degree $-1$ and degree $1$ generators of $\pz{V}$ respectively. We define  comultiplication by
 \bq \Delta(\pz{X})=\pz{X}\otimes \pz{X}\qquad
 \Delta(\pz{1})=\pz{1}\otimes \pz{X}+\pz{X}\otimes \pz{1}.\eq
 The unit map $\iota:\bb{Z}\rightarrow \pz{V}$ by $\iota(1)=\pz{1}$. The trace map is defined by
 \bq\epsilon(\pz{X})=1\qquad
  \epsilon(\pz{1})=0.\eq
It is evident that multiplication is given by \bq m(\pz{1},\pz{X})=m(\pz{X},\pz{1})=\pz{X}\qquad m(\pz{X},\pz{X})=0\qquad m(\pz{1},\pz{1})=\pz{1}.\eq
 Definitions above are made by choosing a basis for $H^*(S^2)$. In section \ref{Hnmodulestr} we give a definition which does not need the choice of a basis.
The TQFT $\mathcal{F}$ assigns to each closed one dimensional manifold (i.e. a circle)  the vector space $\pz{V}$, to each cap the unit map $\iota:\bb{Z}\rightarrow \pz{V}$, to each cup the trace map $\epsilon: \pz{V}\rightarrow \bb{Z}$, to each pair of pants the  multiplication $m:\pz{V}\otimes\pz{V}\rightarrow \pz{V}$,  
and to each reverse pair of pants the comultiplication $\Delta:\pz{V}\rightarrow\pz{V}\otimes \pz{V}$.

\subsection{Tangle cobordisms and the rings $H^m$}\label{ringshn}
We denote the Cartesian coordinates on $\bb{C}\times [0,1]\times[0,1]$ by $(z,t,s)$. 
For a subset $A\subset (\bb{C}\times [0,1]\times[0,1])$ we set $$\partial^v_i A=A\cap (\bb{C}\times [0,1]\times\{i\})$$ and
 $$\partial^h_i A=A\cap (\bb{C}\times\{i\}\times [0,1]).$$
\bdf Let $T_0, T_1$  be two $(m,n)$-tangles. A cobordism between $T_0$ and $T_1$ is a smoothly embedded surface $S$ in $\bb{C}\times [0,1]\times[0,1]$ s.t.
$$\partial^v_i S=T_i$$ for $i=0,1.$
 We also require $S$ to be the product of $\partial^h_i$ or $\partial^v_i $ with a small subinterval  in a \nbhd of each face of $\bb{C}\times \partial([0,1]\times[0,1])$.
\edf
The identity cobordism between $T$ and itself is denoted by $\mathbf{1}_T$.
%
%
Tangle cobordisms can be composed in two ways. First the \textit{vertical composition}: if $S,S'$ are cobordisms between $T_0,T_1$ and $T_1,T_2$ then we get a cobordism
 \bq S'\circ S=\frac{S' \cup S} {\partial^v_0S'\sim \partial^v_1 S} 
\eq between $T_0$ and $T_1$. Secondly the \textit{horizontal composition}: if $S$ is a cobordism between $k T_0 l$ and   $mT_1n$, and $S'$  is a cobordism between $l T'_0 L$ and  \:$n T'_1 N$
 then we get a cobordism 
 \bq
S'S= \frac{S' \cup S} {\partial^h_0S'\sim \partial^h_1 S} 
\eq
 between $kT'_0\circ T_0L$ and $mT'_1\circ T_1N$. The last assumption in the definition of a cobordism ensures that compositions are smooth embedded surfaces.


%

\begin{figure}[ht]
 \includegraphics[scale=.4, bb=0 0 784 172]{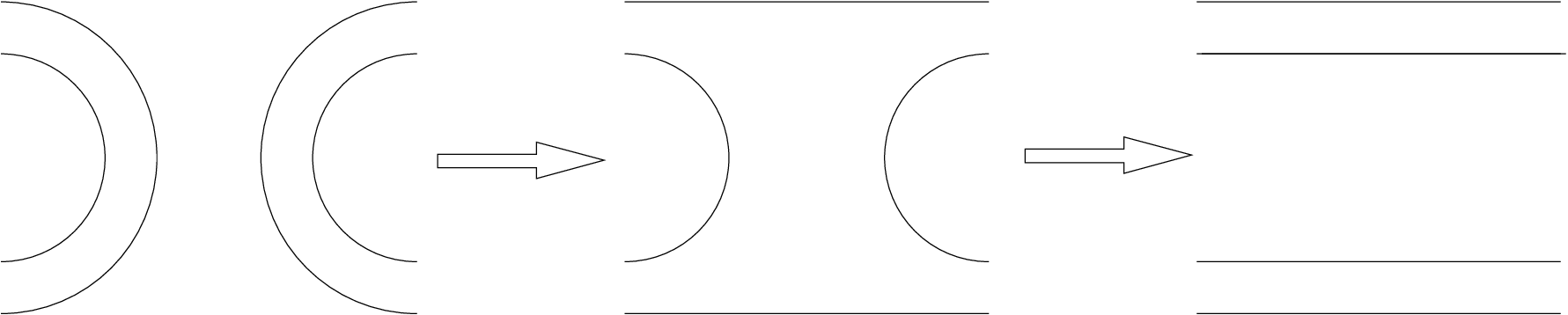}
 \caption{A minimal cobordism}
 \label{minimal_cobo}
\end{figure}

For the purpose of this  paper we just need to consider a special class of tangle cobordisms.
\bdf\label{defminimalcobo}
For a crossingless matching $a\in \P_m$, the \textit{minimal cobordism}  between $a^ta$ and $id_m$ is the one which is given by merging the corresponding strands of $a^t$ and $a$ from the outermost one to the innermost one as depicted in Figure \ref{minimal_cobo}. We denote this minimal cobordism by $S_a$.
\edf

For $a,b\in\pz{C}_m$ set \bq\label{aHb}_a H^m_b=\mathcal{F}(a^tb)\{m\}\eq  and \bq H^m=\bigoplus_{a,b\in\pz{C}_m}  \phantom{a}_a H^m_b.\eq Note that $a^tb$ is a disjoint union of circles so $_a H^m_b=\pz{V}^{\otimes k}$ where $k$ is the number of the circles. Multiplication $$_a H^m_b\otimes \,_cH^m_d\rightarrow \,_aH^m_d$$ is defined to be zero if $b\neq c$. If $b=c$, let $S_b$ be the minimal cobordism between $b^tb$ and $id_m$. The cobordism $\mathbf{1}_a^t S_b \mathbf{1}_c$ is a surface without corners so we get a map $\mathcal{F}(\mathbf{1}_a^t S_b \mathbf{1}_c):\, _a H^m_b\otimes \,_b H^m_c\rightarrow \,_aH^m_c$ which gives us multiplication.\\


We now recall a recursive decomposition of $H^m$ from \cite{RR1}.
Denote by $\P_m'$ the subset of $\P_m$ consisting of elements which contain $\cap_1$, i.e. elements which contain an arc between points 1 and 2,  and denote by $\P_m''$ its complement. (1 can be replaced with any $1\leq i\leq 2m-1$.) $\P_m'$ is in one-to-one correspondence with $\P_{m-1}$. We have a map $\P_m''\rightarrow \P_{m-1}$, $a\mapsto a'$, given by joining the two strands of $a$ that stem from 1 and 2. Let $\P\P_m^1\subset \P_m''\times \P_m''$ be the subset of all $(a,b)$ such that the arcs passing through points $1$ and $2$ in $a^tb$ form a single circle. If $(a,b)$ is in the complement of $\P\P_m^1$ then the arcs passing through points $1$ and $2$ in $a^tb$ form two circles.  

Let  $a,b\in \P_m'$. If $\bar{a}$ denotes $a$ with the $\cap_1$ removed then we have 
$$ _aH^m_b=\phantom{H}_{\bar{a}}H^{m-1}_{\bar{b}}\otimes \pz{V}\{1\}.$$
This contributes a summand of $H^{m-1}\otimes\pz{V}\{1\} $ to $H^m$.
Set 
$${\tilde{H}^{m}}=\bigoplus_{a\in \P_m', b\in\P_m''} \mathcal{F}(a^t,b)\{m\}.$$  

The embedded circle $C$ in $a^t b$  which passes through points 1 and 2 contributes a factor of $\pz{V}\{1+i\}$ to
${\tilde{H}^{m}}$ 
 where $i$ is the number of other pairs of points $2k-1,2k$ which $C$ passes through. 
We can set
 $$\tilde{H}^{m}=\bar{H}^m\otimes \pz{V}\{1\}$$
where $\pz{V}\{1\}$ is the ``local'' contribution of the circle containing $\cap_1$ or $\cup_1$. This means that if $a\in \P_m', b\in\P_m''$ and we modify the strands of $a^t b$ passing through 1 and 2 only in a small \nbhd of the points 1 and 2 then we alter only the second factor in $ \bar{H}^m\otimes \pz{V}\{1\} $.
Also denote by $\tilde{H}^m_1$ and $\tilde{H}^m_2$ the contribution of $\P\P_m^1$ and its complement to $H^m$. Again we can write $\tilde{H}^m_1=H^m_1\otimes \pz{V}\{1\}$ and $\tilde{H}^m_2=H^m_2\otimes \pz{V}\{1\}\otimes \pz{V}\{1\}$ where $\pz{V}\{1\}$ resp. $\pz{V}\{1\}\otimes\pz{V}\{1\}$ are ``local''contributions from the single circle resp. the two circles formed by arcs passing through 1 and 2.
Therefore we get

\bq \label{HmHm-1} H^m=\Bigl(\left(H^{m-1}   \oplus \bar{H}^m \oplus \bar{H}^m\oplus H^m_1 \right)\otimes\pz{V}\{1\} \Bigr)\bigoplus H^m_2\otimes\pz{V}\{1\}\otimes \pz{V}\{1\}.\eq
as abelian groups.

\subsection{The Khovanov invariant for flat tangles}\label{khflattang}
\bdf A flat tangle is a tangle which can be embedded into the plane i.e. a tangle without crossings.
\edf

For a flat $(m,n)$ tangle $T$ we define
\bq\kh(T)=\bigoplus_{a\in\P_m,b\in\P_n} \mathcal{F}(a^t T b)\{n\}.\eq

Obviously $H^m=\kh(id_m)$ as abelian groups. The abelian group  $\kh(T)$ for an $(m,n)$-tangle $T$ has the structure of a $(H^m,H^n)$-bimodule which is given by
\bq \mathcal{F}(\mbf{1}_a S_b \mbf{1}_T\mbf{1}_c):  \,_aH^m_b\otimes \mathcal{F}(b^t T c) \rightarrow \mathcal{F}(a^t T c) \eq  for each $a,b,c$ and zero map
 $_aH^m_b\otimes \mathcal{F}(c^t T d) \rightarrow \mathcal{F}(a^t T d) $ if $b\neq c$.

 For the unlinked union $T\sqcup S^1$ we have \bq \kh(T\sqcup S^1)=\kh(T)\otimes \pz{V}=\kh(T)\{1\}\oplus\kh(T)\{-1\} \eq
 If $S$ is a cobordism between two flat $(m,n)$-tangles $T_0,T_1$, it induces a bimodule map \bq\kh(S):\kh(T_0)\rightarrow \kh(T_1)\eq which is given on each component by
 \bq\mathcal{F}(\mbf{1}_a^t S \mbf{1}_b):\mathcal{F}(a^t T_0 b)\rightarrow \mathcal{F}(a^t T_1 b). \eq
The fact that $\kh(S)$ is independent of the isotopy class of $S$ follows from the same property for the TQFT $\mathcal{F}$.
%
 It is obvious that
 \bqa \kh(S\circ S')&=&\kh(S)\circ\kh(S')\eqa

\bl [\cite{functorvalued}] \label{tensorforflat}If $lTm$ and $mT'n$ are flat tangles then $$\kh(T_0\circ T_1)=\kh(T_0)\otimes_{H^m} \kh(T_1)$$ as $(H^l,H^n)$-bimodules.
I f $T_0S_0T_0'$ and $T_1S_1T_1'$ are minimal cobordisms then \bq \kh(S_0S_1)=\kh(S_0)\otimes_{H^m} \kh(S_1).\eq
\el



\subsection{The Khovanov invariant in general}\label{khgeneral}
In this section we present  Khovanov's invariant for general tangles in a roundabout way which is shorter and suitable for our purpose. For a general tangle $T$, $\kh(T)$ is a chain complex of graded bimodules over the rings $H^m$ so it is doubly graded. For a flat tangle $T$  the chain complex
 \bq\cdots \rightarrow 0\rightarrow\kh(T)\rightarrow0\rightarrow\cdots\eq with $\kh(T)$ in (first or homological ) degree  zero. We denote upward shift in first  by $\{i\}$ and downward shift in second grading by $[i]$.
The only elementary braids which are not flat are the braids $\sig^+_{i;m}$ and $\sig^-_{i;m}$. Consider the chain complexes

 \bq C^+_{i;m}\qquad \qquad \cdots\rightarrow 0\rightarrow\kh(id_m) \stackrel{\alpha}\rightarrow \kh(\cup_{i;m}\cap_{i;m}) \{-1\}\rightarrow 0\rightarrow\cdots\eq
  \bq C^-_{i;m}\qquad \qquad \cdots\rightarrow 0\rightarrow\kh(\cup_{i;m}\cap_{i;m})\stackrel{\beta}\rightarrow \kh(id_m)\{-1\}\rightarrow 0\rightarrow\cdots\eq
%
%
 where the domain of maps $\alpha$ and $\beta$ are in (first or homological) degree zero. The map $\alpha$ is $\kh(S_i)=\oplus_{a,b\in\P_m} \mathcal{F}(\mbf{1}_a^t S_i \mbf{1}_b)$ where $S_i$ is the minimal cobordism between $\cup_{i;m}\cap_{i;m}$ and $id_m$. The map $\beta$ is obtained in the same way from $S^t_i$ which is a cobordism between $id_m$ and $\cup_{i;m}\cap_{i;m}$. The $-1$ degree shift is to make the map $\alpha$ of (the second or quantum) degree zero.

  Let $\sig^+=\sig^+_{i;m}$ and $\sig^-=\sig^-_{i;m}$ be as in the Figure \ref{posnegbraid}. 
  Khovanov defines
 \bqa\label{sigchain} \kh(\sig^+)&=& C^+_{i;m}\{-1\}\\
 \label{sigchain2}\kh(\sig^-)&=& C^-_{i;m}[1]\{2\}.
 \eqa

 Now let $n_0T_0n_1T_1n_2\cdots n_kT_kn_{k+1}$ be a decomposition of a tangle $T$ into elementary tangles.
 \bdf \bq\label{khdef}\kh(T):=\kh(T_0)\otimes_{H^{n_1}}\kh(T_1)\otimes_{H^{n_2}}\cdots \otimes_{H^{n_k}}\kh(T_k)  \eq
 \edf

In \cite{functorvalued}, Khovanov defines his invariant using the cube of resolutions and obtains the above equation as a consequence. He also shows that $\kh(T)$ is independent of the decomposition and is invariant under isotopies of $T$.
If $L$ is a link, the homology  of $\kh(L)$ is the original Khovanov homology \cite{categorification}   of $L$ with its first grading reversed. We set
\bq \overline{\kh}^i(T)=\bigoplus_{j-k=i} H(\kh(T))^{j,k}.\eq
Seidel and Smith \cite{SS} conjecture that their invariant $\khs$ is equivalent to $\overline{\kh}$.

\bpr\label{khofbraid}
\begin{enumerate}
\item With the same notation as in \eqref{HmHm-1} we have $$\kh(\cup_{i;m}\cap_{i;m} )=\left(\left(H^{m-1}\oplus \bar{H}^m\oplus \bar{H}^m \oplus H^m_1\right)  \otimes \pz{V}\{1\}\otimes \pz{V}\right)\bigoplus H^m_2\otimes\pz{V}\{2\}.$$  
\item On the first four direct summands of \eqref{HmHm-1}, the map $\alpha:  \mathcal{F}(a^t id_m b) \rightarrow \mathcal{F}(a^t \cup_{i;m}\cap_{i;m} b)$ is given by the comultiplication $\Delta:\pz{V}\{1\} \rightarrow \pz{V}\{1\}\otimes\pz{V}$ tensored with the identity map. On the last one it is given by the multiplication $m:\pz{V}\{1\}\otimes\pz{V}\{1\}\ra\pz{V}$ tensored with the identity.  
\item On the first four direct summands, the map $\beta: \mathcal{F}(a^t \cup_{i;m}\cap_{i;m} b)\rightarrow \mathcal{F}(a^t id_m b)$ is given by the multiplication $m:\pz{V}\otimes\pz{V}\rightarrow\pz{V}$ tensored with the identity map. On the last one it is given by the comultiplication $\Delta$ tensored with the identity.  

\item\bqas
\overline{\kh}(\sigma^+_k)&=&  \left(   H^{m-1}\oplus \bar{H}^m\oplus \bar{H}^m  \oplus H^m_1\right)\otimes \pz{V}\{-1\}   \bigoplus H_2^m \otimes \pz{V}\{2\}\\
\overline{\kh}(\sigma^-_k)&=& \left(H^{m-1}\oplus \bar{H}^m \oplus \bar{H}^m  \oplus H^m_1 \right)\otimes \pz{V}\{3\} \bigoplus H^m_2\otimes\pz{V}\{2\}
\eqas
\end{enumerate}
\epr
\bp
 i) Follows easily by comparison to \eqref{HmHm-1}.\\
ii) This is because in the first four summands the cobordism merges two circles into one and in the last one it decomposes a circle into two.\\
iii) Similarly because in the first four summands the cobordism decomposes one circle into two and in the last one it merges two circles into one.\\
iv) Since $m:\pz{V}\otimes\pz{V}\rightarrow\pz{V}$ is surjective we have
$$H^1(C^-_{i;m})= \frac{H^m_2\otimes\pz{V}\{1\}\otimes \pz{V}\{1\}\{-1\}}{ H^m_2\otimes \on{Im}\Delta}=H^m_2\otimes \frac{\pz{V}\otimes \pz{V}\{1\}}{\on{Im}\Delta}.$$
$\pz{V}\otimes \pz{V}\{1\}/\on{Im}\Delta$ is isomorphic to $\bb{Z}<\pz{1}\otimes \pz{1}, \pz{1}\otimes \pz{X}-\pz{X}\otimes \pz{1}>\{1\}$.
   Therefore $$H^1(C^-_{i;m})= H^m_2\{-1\}\oplus H^m_2\{1\}\cong H^m_2\otimes\pz{V}.$$

   The map $\Delta$ is injective and the kernel of
   $m:\pz{V}\{1\}\otimes\pz{V}\ra\pz{V}\{1\}\{-1\}$ equals $\Z \langle \pz{X}\otimes\pz{X},\pz{1}\otimes \pz{X}-\pz{X}\otimes \pz{1} \rangle$. Tensoring with $\pz{X}\otimes\pz{X}$ has the effect of shifting degree by two, and tensoring with $\pz{1}\otimes \pz{X}-\pz{X}\otimes \pz{1}$ does not shift the degree.  Therefore
   $$H^0(C^-_{i;m})=\left(H^{m-1}\{1\}\bigoplus \bar{H}^m \{1\}\bigoplus \bar{H}^m \{1\}\bigoplus H^m_1\{1\}\right)\otimes  (\Z\{2\}\oplus \Z )$$ 
   which is isomorphic to
   $ \left(H^{m-1}\bigoplus \bar{H}^m \bigoplus \bar{H}^m  \bigoplus H^m_1 \right)\otimes \pz{V}\{2\}.$

Therefore
\bqas
 \overline{\kh}(\sigma^-_k) =&(H^0(C^-_{k;m})\oplus H^1(C_{k;m})\{1\})\{+1\}\\
 =& \left(H^{m-1}\oplus \bar{H}^m \oplus \bar{H}^m  \oplus H^m_1 \right)\otimes \pz{V}\{2+1\} \bigoplus H^m_2\otimes\pz{V}\{2\}
 \eqas
Now for $\sig_i^+$ we have
\bqas
 H^0(C_{i;m}^+)=& H_2^m\otimes(\ker m: \pz{V}\{1\}\otimes \pz{V}\{1\} \rightarrow \pz{V}\{2\}\{-1\}) \\
 =&H_2^m\otimes \Z< \pz{X}\otimes\pz{X},\pz{1}\otimes \pz{X}-\pz{X}\otimes \pz{1}> \{ 2\}  \\ \cong & H_2^m \otimes \pz{V}\{3\}.
\eqas
\bqas
  H^1(C_{i;m}^+)= &\left(   H^{m-1}\bigoplus \bar{H}^m\bigoplus \bar{H}^m  \bigoplus H^m_1\right) \otimes \coker(\Delta: \pz{V}\{1\}\rightarrow \pz{V}\{1\}\otimes \pz{V} )\{-1\} \\
  = &\left(   H^{m-1}\bigoplus \bar{H}^m\bigoplus \bar{H}^m  \bigoplus H^m_1\right)\otimes \{ \pz{1}\otimes \pz{1}, \pz{1}\otimes \pz{X}-\pz{X}\otimes \pz{1} \} \\
   \cong &  \left(   H^{m-1}\bigoplus \bar{H}^m\bigoplus \bar{H}^m  \bigoplus H^m_1\right)\otimes \pz{V} \{-1\}
\eqas
\bqas
\overline{\kh}(\sigma^+_k)= & (H^0(C^+_{k;m})\oplus H^1(C_{k;m})\{1\})\{-1\}\\
= & \left(   H^{m-1} \oplus \bar{H}^m \oplus \bar{H}^m  \oplus H^m_1\right)\otimes \pz{V}\{-1\}
 \bigoplus H_2^m \otimes \pz{V}\{3-1\}
\eqas

\ep

\section{The symplectic Khovanov homology of links}\label{SSRRChapter}

  In this  section we review the construction of Seidel and Smith \cite{SS}.  
This material appeared  in \cite{RR1} as well.
Denote by $\mathrm{Conf}_m$ the space of all unordered $m$-tuples of  distinct complex numbers   $(z_1,\cdots,z_m).$ Denote by $\mathrm{Conf}^0_m$ the subset of  $\mathrm{Conf}_m$ consisting of $m$-tuples which add up to zero, i.e. $z_1+\cdots+z_m=0$.
\renewcommand{\o}{\mathbf{0}_{2\times 2}}
\renewcommand{\i}{\mathbf{1}_{2\times 2}}                                         
Let  $\i$ and $\o$ denote the $2\times 2$ identity matrices respectively. 
 Let $\mathcal{S}_m$ be the set of matrices in $\mathfrak{sl}_{2m}$ of the form  \\ 
\bq\label{matrixS}\left(
\begin{array}{cccccc}
                                           y_1 & \i &\space & & \space& \\
                                           y_2 & & \i &  & & \\
                                           \vdots&  &  & \ddots & \\
                                             y_{m-1}& & & & \i\\
                                             y_m &  &  & & \o\\
                                         \end{array}
                                       \right)\eq\\
 Here $y_1\in \mathfrak{sl}_2$ and $y_i\in \mathfrak{gl}_2$ for $i>1$. 
The set $\Sm$ is in fact a homogeneous slice transverse to the orbit of $x$ (\cite{SS}, Lemma 23).
The map  $\chi$ restricted to $\mathrm{Conf}^0_{2m}$ 
  is a differentiable fiber bundle (\cite{SS}, Lemma 20). 
   We denote the fiber of $\chi$ over $t$ by $\mathcal{Y}_{m,t}$, i.e. $$\mathcal{Y}_{m,t}=\chi^{-1}(t)$$ If $t=(\mu_1,\ldots,\mu_{2m}) \notin \mathrm{Conf}^0$, by  $\mathcal{Y}_{m,t}$ we mean $\mathcal{Y}_{m,t'}$  where $t'=(\mu_1-\sum \mu_i/2m,\ldots, \mu_{2m}-\sum \mu_i/2m).$  Each such fiber inherits a  a K\"ahler  structure from $\Sm$ and this equips the total space with a connection in a canonical way. 
Since these fibres are not compact the existence of parallel transport maps is not guarantied. By modifying the Ka\"hler structure on $\Sm$ Seidel and Smith assign rescaled parallel transport maps
\bq\label{paralleltransD}h_\beta^{res}:\mathcal{Y}_{m,\beta(0)}\ra\mathcal{Y}_{m,\beta(1)}\eq
 to each curve $\beta:[0,1]\ra\mathrm{Conf}_{2m}$.
The map $h^{res}_{\beta}$  is a symplectomorphism defined on arbitrarily large compact subsets of $\mathcal{Y}_{m,\beta(0)}$.

Let $E_{y}^\mu$ denote the $\mu$-eigenspace of $y$.
\begin{lemma}[\cite{SS}, Lemmata 25 and 26 ]\label{2dkernel}
  For any $y\in \Sm$  and $\mu\in \mathbb{C}$ the projection $\mathbb{C}^{2m}\rightarrow \mathbb{C}^2$ onto the first two coordinates gives an injective map $E_y^\mu\rightarrow \mathbb{C}^2$.
  Any eigenspace of any element $y\in \Sm$ has dimension at most two. Moreover the set of elements of $\mathcal{S}_m$ with 2 dimensional kernel can be canonically identified with $\mathcal{S}_{m-1}$ and this identification is compatible with $\chi$.
\end{lemma}

For a subset $A\subset \mathfrak{sl}_{2m} $, let  $A^{sub,\lambda}$ (resp. $A^{sub3,\lambda}$) be the subset of matrices in $A$ having eigenvalue $\lambda$  of multiplicity two (resp. three) and two Jordan blocks of size one (resp. two Jordan blocks of sizes 1,2) corresponding to the eigenvalue $\lambda$  and no other coincidences between the eigenvalues. Here are two results describing  neighborhoods of $\mathcal{S}_m^{sub,\lambda} $ and  $\mathcal{S}_m^{sub3,\lambda} $ in $\Sm$.

\begin{lemma}[\cite{SS}, Lemma 27]\label{2nbhd}
Let $D\subset \on{\mathrm{Conf}}^0_{2m}$ be a disc consisting of the $2m$-tuples $(\mu-\varepsilon, \mu-\varepsilon,\mu_3,\ldots, \mu_{2m})$ with $\varepsilon \in \C$ of small magnitude. Then there is a neighborhood $U_\mu$ of $\mathcal{S}_m^{sub,\mu} $ in $\Sm\cap \chi^{-1}(D)$ and  an isomorphism $\phi$ of $U_\mu$  with  a \nbhd of $\Sm^{sub,\mu}$ in $\Sm^{sub,\mu}\times \mathbb{C}^3$ such that $f \circ \phi=\chi$ on $\Sm\cap \chi^{-1}(D)$ where $f(x,a,b,c)=a^2+b^2+c^2$.
  Also if $N_y \Sm^{sub,\mu}$ denotes the normal bundle to $\Sm^{sub,\mu}$ at $y$ then we have \bq\label{2nbhdmap}\phi(N_y \Sm^{sub,\mu})= \mathfrak{sl}(E_{y}^{\mu})\oplus \zeta_y\eq where 
 $\zeta_y$ is the trace free part of $\{\mathbb{C}\cdot 1\subset \mathfrak{gl}(E_y^\mu) \} \oplus \mathfrak{gl}(E_y^{\mu_3})\oplus \ldots \oplus \mathfrak{gl}(E_y^{\mu_{2m}}) $.   
\end{lemma}

  Consider the line bundle $\mathcal{F}$ on $\Sm^{sub3,\mu}$ whose fiber at $y\in\Sm^{sub3,\mu}$ is $(y-m)E_{y_s}^\mu$ where $y_s$ is the semisimple part of $y$. To $\mathcal{F}$ one  associates a  $\Bbb{C}^4$ bundle $\mathcal{L}=(\mathcal{F}\backslash 0)\times_{\mathbb{C}^*} \mathbb{C}^4$ where  $z\in \mathbb{C}^*$ acts on $\mathbb{C}^4$  by
  \begin{equation}\label{C2action}
  (a,b,c,d) \rightarrow (a,z^{-2}b, z^{2}c, d ).
\end{equation}
$\mathcal{L}$ decomposes as \bq\label{Ldecomp}\mathcal{L}\cong\mathbb{C}\oplus\mathcal{F}^{-2}\oplus\mathcal{F}^2\oplus \mathbb{C}.\eq
Fibers of $\mathcal{L}$ should be regarded as transverse slices in $\frak{sl}_3$. Upon choosing  suitable coordinates on such a transverse  slice (at the zero matrix) the function $\chi$  equals the function $p:\frak{sl}_3\rightarrow \mathbb{C}^2$ given by \bq\label{p} p(a,b,c,d)=(d,a^3-ad+bc).\eq  $p$ is also well-defined as a function $\mathcal{L}\rightarrow \mathbb{C}^2$ because $b$ and $c$ are coordinates on line bundles which are inverses of each other.  
Denote by $\tau(d,z)$ the set of solutions of  $\lambda^3-d\lambda+z=0$.

\begin{lemma}[\cite{SS}, Lemma 28]\label{3nbhd}
 Let $P\subset \mathrm{Conf}^0_{2m}$ be the set of $2m+2$-tuples \bq\label{3nbhddisk}(\mu_1,\ldots, \mu_{i-1},\tau(d,z),\mu_{i+3},\ldots ,\mu_{2m+2}).\eq where $d$ and $z$ vary in a small disc in $\mathbb{C}$ containing the origin.
 There is a \nbhd  $V$ of $\Sm^{sub3}$ in $\Sm\cap\chi^{-1}(P)$ and an isomorphism $\phi'$ from $V$ to a \nbhd of zero section in $\mathcal{L}$ such that  $p(\phi'(x))=(d,z)$ if \\
 $$\chi(x)=(\mu_1,\ldots, \mu_{i-1},\tau(d,z),\mu_{i+3}, \ldots, \mu_{2m+2}).$$
\end{lemma}

If $\mu\in \mathbb{C}^{2m}/\mathbf{S}_{2m}$ has only one element of multiplicity two or higher, which we denote by $\mu_1$, denote by $\mathcal{D}_{m,\mu}$ the set  of singular elements  of $(\chi^{-1}(\mu)\cap \mathcal{S}_m)$ i.e.
\bq\label{Ddef}\mathcal{D}_{m,\mu}=(\chi^{-1}(\mu)\cap \mathcal{S}_m)^{sub,\mu_1}.\eq
 Let $\mathcal{D}_m$ be the union of all these $\mathcal{D}_{m,t}$ regarded as a subset of $\mathcal{S}_m$. It inherits a K\"ahler metric from $\mathcal{S}_m$.  We have the map $\chi: \mathcal{D}_{m}\rightarrow \mathbb{C}\times \mathbb{C}^{2m-2}/\mathbf{S}_{2m-2} $. By forgetting the first eigenvalue, the image of $\chi$ can be identified with $\mathrm{Conf}_{2m-2}$.

\subsection{Relative vanishing cycles }\label{relvanish}

Let $X$ be a complex manifold and $K$ a compact submanifold. Let  $g$ be a K\"ahler metric on $Y=X \times \mathbb{C}^3$ (not necessarily the product metric)  and denote its imaginary part by $\Omega$. Consider the map $f:X\times \mathbb{C}^3\rightarrow \mathbb{C}$ given by $f(x,a,b,c)= a^2+b^2+c^2$ and denote by $\phi_t$ the gradient flow of $-\on{Re} f$. 
Let $W$ be the set of points  $y\in Y$ for which the trajectory $\phi_t(y)$ 
 exists for all positive $t$. 
One can prove that 
$W$ is a manifold and the mapping $l:W\rightarrow X$ given by $l(y)=\lim_{t\rightarrow \infty} \phi_t(y)$ is well-defined and smooth.   We have $\Omega|_W=l^* \Omega|_X$. The function $f$ restricted to $W$ is real and nonnegative.

Set $V_t(K)=\pi^{-1}(t)\cap l^{-1}(K)=l|_{\pi^{-1}(t)}^{-1}(K)$ which is a manifold for $t$ small.  
It follows from Morse-Bott lemma that $V_t(K)$ is a 2-sphere bundle on $K$ for $t$ small. To generalize the invariant to tangles we will need a slightly more general version of the above construction in which $K$ is noncompact and the metric equals the product metric outside a compact subset. (See subsection \ref{func}.) The resulting vanishing cycle equals (symplectically) the product bundle outside a compact subset.

\subsection{Fibered $A_2$  singularities}\label{fibered}

Assume we have the same situation as in the Lemma \ref{3nbhd}, i.e. let $\mathcal{F}$ be a holomorphic line bundle over a complex manifold $X$ and define $Y$ to be $(\mathcal{F}\backslash0) \times_{\mathbb{C}^*} \mathbb{C}^4$ where the $\mathbb{C}^*$ action is as in the formula \eqref{C2action}. Let $\Omega$ be an arbitrary K\"ahler form on $Y$ and by  regarding $X$ as the zero section of $Y$, $\Omega$ restricts to a K\"ahler form on $X$. Let $(a,b,c,d)$ be the coordinates on fibers of $Y\rightarrow X$ and $(d,z)$ coordinates on $\mathbb{C}^2$. Let the map $p:Y\rightarrow \mathbb{C}^2$ be as in Lemma \ref{3nbhd}. Let $Y_d=p^{-1}(\mathbb{C}\times \{d\})$ and $p_d: Y_d\rightarrow \mathbb{C}$ be the restriction of $p$. Set $Y_{d,z}=p^{-1}(d,z)$. For $d\neq 0$ the critical values of $p_d$ are \smash{$\zeta^{\pm}_d=\pm 2\sqrt{d^3/27}$}.

 Let $K$ be Lagrangian submanifold of $X$.  Using the relative vanishing cycle construction for the function $p_d$ we can obtain a Lagrangian submanifold $L_d$ of $Y$ which is a  sphere bundle over $K$. (This construction works when $Y$ is a nontrivial bundle over $X$ as well.) There is another way of describing this  Lagrangian as follows. Let $\mathrm{Y}\cong \mathbb{C}^4$  be the fiber of $Y\rightarrow X$ over some point of $X$ and let $p:\mathrm{Y}\rightarrow \mathbb{C}^2$ be as before. The restriction of the $\mathbb{C}^*$ action to $S^1$ is a Hamiltonian action with the moment map $\mu(a,b,c,d)=|c|^2-|b|^2$.
 Define
 \bq\label{C} C_{d,z,a}=\{(b,c):  \mu(b,c)=0, a^3-da-z=-bc \}\subset \mathrm{Y}_{d,z} \eq which is a point if $a^3-da-z= 0$ and a circle otherwise.  The three solutions of this equation correspond to the critical values of the projection $q_{d,z}:\mathrm{Y}_{d,z}\rightarrow \mathbb{C}$ to the $a$ plane. In the situation of Lemma \ref{3nbhd} they correspond to the three  eigenvalues of $\mathrm{Y}$. Let $\alpha(r)$ be any embedded curve in $\mathbb{C}$ which intersects these critical values (only) if $r=0,1$. Define \bq\label{Lambd}\Lambda_\alpha=\bigcup_{r=0}^{1}  C_{d,z,\alpha(r)}\eq which is a Lagrangian submanifold of $\mathrm{Y}_{d,z}$ (with K\"ahler form induced from $\mathbb{C}^4$). Let $\mathbf{c}, \mathbf{c}',\mathbf{c}''$ be as in the Figure \ref{curvesalphabeta} where dots represent the critical values of $q_{d,z}$.
We can associate to $K$ a Lagrangian submanifold $\Lambda_{d,\alpha}$ of $Y$ by defining $\Lambda_{d,\alpha}=(Y|K)\times_{S^1} \Lambda_{\mathbf{c}}$.  
  Seidel and Smith prove that these two procedures give the same result (\cite{SS}, Lemma 40):

 \begin{lemma}\label{2lagrang}
 If the K\"ahler form on $Y$ is obtained from  a K\"ahler form on $X$, a Hermitian metric on $\mathcal{F}$ and the standard form on $\mathbb{C}^4$ then $L_d=\Lambda_{d,\mathbf{c}}$.

 \end{lemma}

\begin{figure}[ht]
\begin{center}
\scalebox{.7}{
\input{curvesalphabeta.pstex_t}
}
\caption{}
 \label{curvesalphabeta}
\end{center}
\end{figure}

\subsection{Lagrangian submanifolds from crossingless matchings}\label{crossless}

\begin{figure}[ht]
\includegraphics[width=4.5in]{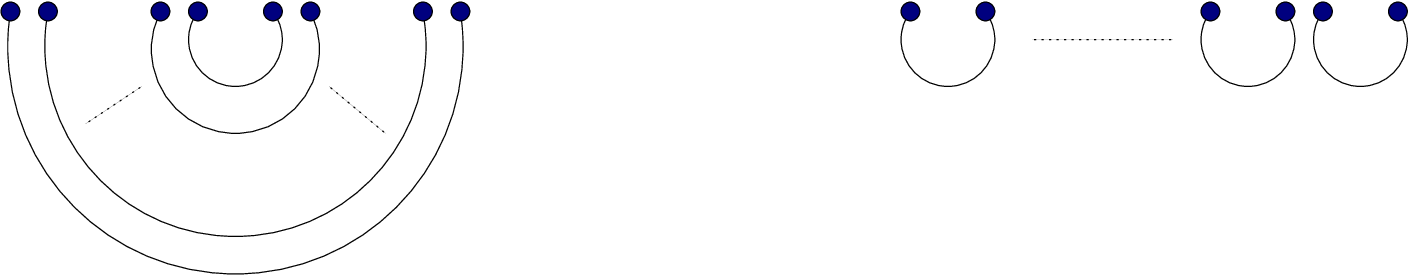}
\caption{Two crossingless matchings}
\label{crossinglessmatchings}
\end{figure}

Let $\mu\in \mathrm{Conf}_{2m}$. A \textit{crossingless matching} $D$ with endpoints in $\mu$ is a set of $m$ disjoint embedded curves $\delta_1,\ldots,\delta_m$ in $\mathbb{C}$ which have (only) elements of $\mu$ as endpoints.  See Figure \ref{crossinglessmatchings}. To $D$ we associate a Lagrangian submanifold $L_D$ of $\mathcal{Y}_{m,\mu}$ as follows. Let 
$\{\mu_{2k-1},\mu_{2k}\}\subset \mu$  be the endpoints of $\delta_k$ for each $k$. Let $\gamma$ be a curve in $\mathrm{Conf}^0_{m}$ such that  $\gamma(t)=\{\gamma_1(t), \gamma_2(t), \mu_3, \mu_4,\ldots,\mu_{2m}\}$, $\gamma_i(0)=\mu_i, i=1,2$ and as $s\rightarrow 1$, $\gamma_1(t), \gamma_2(t)$ approach each other on $\delta_1$  and collide. For example if $\delta_1(t)$ is a parameterization of $\delta_1$ s.t $\delta_1(0)=\mu_1,\delta_1(1)=\mu_2$ the we can take $\gamma(t)=\{\delta(t/2), \delta(1-t/2), \mu_3,\ldots,\mu_{2m}\}$. Set $\bar{\mu}=\mu \backslash \{\mu_{1},\mu_{2}  \}$, $\mu'=\gamma(1)$.

If $m=1$ then relative vanishing cycle construction for $\chi: \mathcal{S}_1\rightarrow \bb{C}$ with the critical point over $\gamma(1)=0$ gives us  a Lagrangian sphere $L$ in $\mathcal{Y}_{1,\gamma(1-s)}$ for small $s$. Using reverse parallel transport along $\gamma$ we can move $L$ to $\mathcal{Y}_{1,\mu}$ to get our desired Lagrangian submanifold.  Now for arbitrary $m$ assume by induction that  we have obtained a Lagrangian $L_{\bar{D}}\subset \mathcal{Y}_{m-1,\bar{\mu} } $ for  $\bar{D}$ which is obtained from $D$ by deleting $\delta_1$.
Now $\mathcal{Y}_{m-1,\bar{\mu}}$ can be identified with  $D_{m,\tau}$ where $\tau=(0,0, \mu_3-(\mu_1+\mu_2)/(2m-2),\ldots,  \mu_{2m}-(\mu_1+\mu_2)/(2m-2))$. Use parallel transport  to move $L_{\bar{D}}$ to $D_{m,\gamma(1)}$ . The later one is the set of singular points of $\mathcal{Y}_{m,\gamma(1)}$  so Lemma \ref{2nbhd} tells us that we can use relative vanishing cycle construction for $L_{\bar{D}}$ to obtain a Lagrangian in $\mathcal{Y}_{m,\gamma(1-s)}$ for small $s$. Parallel transporting  it along $\gamma$ back to $\mathcal{Y}_{m,\mu}$ we obtain our desired Lagrangian which is topologically a trivial sphere bundle on $L_{\bar{D}}$. We see that $L_D$ is  diffeomorphic to a product of spheres. Different choices of the curve $\gamma$ result in Hamiltonian isotopic Lagrangians. The same holds if we isotope the curves in $D$ inside $\mathbb{C}\backslash \mu$.\\

Now we can define the Seidel-Smith invariant. 
  Since each manifold $\mathcal{Y}_{m,\nu}$ is a   submanifold of the affine space $\Sm$ and has trivial  normal bundle, its Chern classes are zero. This together with the fact that $H_1(\mathcal{Y}_m)=0$ implies that   the canonical bundle of  $\mathcal{Y}_{m,\nu}$ is trivial and so has a unique infinite Maslov cover. We start by choosing global sections $\eta_{\Sm}$ and $\eta_{\frak{h}/W}$. Then we choose trivializations for regular fibers of $\chi_{\Sm}$  characterized by $\eta_{\mathcal{Y}_{m,t}}\wedge \chi^*\eta_{\frak{h}/W}=\eta_{\Sm}.$ If we choose a grading for $L\subset \mathcal{Y}_{m,t_0}$ and $\beta$ is a curve in $\mathrm{Conf}_{2m}$ starting at $t_0$, one can continue the grading on $L$ uniquely to $h_{\beta|_{[0,s]}}(L)$ for any $s$. Therefore the grading of $L$ uniquely determines that of $h_\beta(L)$. 

Let  $\mathcal{D_+}$ be the crossingless matching at the left hand side of picture \ref{crossinglessmatchings}. If a link $K$ is obtained as closure of a braid $\beta\in Br_m$, let $\beta'\in \mathrm{Conf}_{2m}$ be the braid obtained from $\beta$ by adjoining the identity braid $\on{id}_m$.
\begin{definition}\label{defss}
$$\khs^{*}(K)=HF^{*+m+w}(L_{\mathcal{D_+}}, h^{res}_{\beta'}  (L_{\mathcal{D_+}}))$$
\end{definition}
Here  $w$ is the writhe of the braid presentation, i.e. the number of positive crossings minus the number of the negative crossings in the presentation.
Since the manifold is convex at infinity and the Lagrangians are exact, the above Floer cohomology is well-defined.
Independence from choice of $\beta$ is proved in \cite{SS}, section 5C. 

\section{Generalization to tangles}\label{SSTangle}
In this section we recall the generalization of the invariant of Seidel and Smith to tangles from \cite{RR1}. 

\subsection{The functor associated to a tangle}\label{func}

Let \bq\label{tdecomp}n_0 T_1 n_1 T_2\ldots n_{l-1} T_l n_l\eq be a decomposition of an oriented tangle $T$ into elementary tangles.
%
Set $\nu_j=i(T_j)$ and $\nu_{l+1}=t(T_l)$
 We have $\nu_i\in \mathrm{Conf}_{n_i}$ for $i=0,\ldots, l$.
To each $T_i$ we want to associate a Lagrangian correspondence $L_{i,i+1}=L_{T_i}$ between $\mathcal{Y}_{n_i,\nu_i}$ and $\mathcal{Y}_{n_{i+1},\nu_{i+1}}$. 
In this way we can associate to $T $ a generalized Lagrangian correspondence\\
 \begin{equation}\label{Phi}
 \Phi(T)=(L_{0,1},L_{1,2},\ldots, L_{n-1,n})\{ -m-w\}
 \end{equation}
\\
from $\mathcal{Y}_{n}$ to $\mathcal{Y}_{m}$. Here $m$ and $w$ are the number of cups and the writhe (number of positive crossings minus the number of negative ones) of the decomposition respectively. 

If $T_k$ is an elementary  braid in $Br_{2m}$, we set $L_{T_k}$ to be $ \on{graph} (h^{res}_\beta)$   
regardless of the orientation of the braid. Of course we can extend this definition to any braid.
%
 Let $V_i$ be the relative vanishing cycle  for the map $f$ in Lemma \ref{2nbhd} where $i$th and $(i+1)$th eigenvalues ($\mu_i, \mu_{i+1}$) of $\nu_k$ come together  at some point $\mu$. 
Using  a theorem of T.~Perutz (\cite{PerutzI}, Theorem 2.19) one can describe monodromy maps around singularities of symplectic Morse-Bott fibrations as fibered Dehn twists. 
%
Therefore using the local picture of the Lemma \ref{2nbhd} we see that if we have a subset  $B\subset \mathcal{Y}_m$ for which the naive (non-rescaled) parallel transport map $h_{\sig_i}|_B$ is well-defined then \bq\label{monodromydehn}h_{\sigma_i}\cong\tau_{V_i}.\eq

%

Let ${V_i}_x$ denote the ($S^2$) fiber of $V_i$ over $x$. We grade $\tau_{V_i}$ in such a way that \bq\label{gradedehn}\tau_{V_i} {V_i}_x={V_i}_x\{1\}\eq and the grading function vanishes outside a \nbhd of $V_i$. This grading is unique. (Lemma 5.6 in \cite{gradedlag})

If $T_i=\cup_{j;m}$, we define a  Lagrangian   $L_{\cup_{j;m}}$, regardless of the orientation of $\cup_{k;m}$,  as follows.
The result depends on a real parameter $R>0$.
To simplify the notation we set $k=j, l=j+1$.     
%
 %
With $\nu_i$ as given above let $\nu=\nu_i=\{\mu_1,\ldots, \mu_{2m}\}$.
 Let $\gamma$ be a curve in $\mathrm{Conf}^0_{2m}$ such that $\gamma(0)=\nu_i$ and as $s\rightarrow 1$, $\mu_{k}$ and $\mu_{l}$ approach each other linearly and collide at a point $\mu'$.  For example we can take $$\gamma(t)=\{\mu_1,\ldots,\mu_k+t(\mu_l-\mu_k)/2,\ldots,\mu_l-t (\mu_l-\mu_k)/2 ,\ldots,\mu_{2m} \}$$  provided that $\mu_k+t(\mu_l-\mu_k)/2$ does not intersect the other $\mu_i$. Set $\nu^{k,l}=\nu \backslash \{\mu_{k},\mu_{l}  \}$, $\nu'=\gamma(1)$.  
  We use Lemma \ref{2nbhd} to identify a \nbhd of $\Sm^{sub,\mu'}$ in $\Sm$ locally with $\Sm^{sub,\mu'} \times \mathbb{C}^3$. This induces a K\"ahler form and hence a metric on $\Sm^{sub,\mu'}\times \mathbb{C}^3$. We perturb the complex structure outside a compact ball of radius $\rho$ (to be chosen below) so that outside that set the resulting metric equals the product  metric. Now we use the relative vanishing cycle construction for the whole $\Sm^{sub,\mu'}$.  It  yields (after restriction) a sphere bundle $V=V_{\gamma(1-s)}(\Sm^{sub,\mu'}) \subset \mathcal{Y}_{m,\gamma(1-s)}$ for small $s$  with projection $\pi:V\rightarrow \mathcal{Y}_{m,\nu'}\cap \Sm^{sub,\mu'_k}$. The relative vanishing cycle construction can be used because the metric equals the product metric outside a compact set.

We denote the image of $V$ under  parallel transport map along $-\gamma$, i.e.  
$$  h^{-1}_{\gamma|_{[0,1-s]}}(V)\subset\mathcal{Y}_{m,\nu}$$  by the same notation $V$.
 Composing $\pi$ with the  parallel transport map $h^{-1}_{\gamma|_{[0,1-s]}}$ we obtain a projection $\pi:V\rightarrow  \mathcal{Y}_{m,\nu'}\cap \Sm^{sub,\mu'}$  which is a $S^2$ bundle. By Lemma \ref{2nbhd}, $ \mathcal{Y}_{m,\nu'}\cap \Sm^{sub,\mu'}$ can be identified with $\mathcal{D}_{m-1,\nu'}$  from (\ref{Ddef}). Let $\delta$ be a geodesic in $\mathrm{Conf}^\Sm$  
 joining $\nu'$ to $\nu^{k,l}$.
 We can use parallel transport map \eqref{paralleltransD} along the curve $\delta$ to map  $\mathcal{D}_{m-1,\nu'}$ to $\mathcal{D}_{m-1,\nu^{k,l}\cup \{0,0\}}$. The latter can be identified with $\mathcal{Y}_{m-1,\nu^{k,l}}$. So we obtain a fibration $\pi:V\rightarrow \mathcal{Y}_{m-1,\nu^{k,l}}$. We can use this map $\pi$ to pull $V$ back to $\mathcal{Y}_{m-1,\nu^{k,l}}\times \mathcal{Y}_{m-1,\nu^{k,l}} $. Let ${\cup_{j;m}}$ be its restriction to the diagonal. It  is a Lagrangian submanifold of $\mathcal{Y}_{m,\nu_i}^- \times \mathcal{Y}_{m-1,\nu_i^{k,l}}=\mathcal{Y}_{m,\nu_i}^-\times\mathcal{Y}_{m-1,\nu_{i+1}}$.
 Let $\psi=\psi_1+\psi_2$ be the plurisubharmonic function on $\mathcal{Y}_{m,\nu_i}^-\times\mathcal{Y}_{m-1,\nu_{i+1}}$. We can choose $\rho$ in such a way that the inverse image of $\psi=R$ lies inside the ball of radius $\rho$.
We have a projection $\pi: L_{\cup_{j;m}}\rightarrow \Delta\subset \mathcal{Y}_{m-1,\nu_{i+1}}^-\times\mathcal{Y}_{m-1,\nu_{i+1}}$.\\
  As in the case of Lagrangians from crossingless matchings, 
   replacing the curve $\gamma$  with another curve in the same homotopy class (with fixed endpoints) results in a new $L_{\cup_{j;m}}$ which is Lagrangian isotopic to the former one. Sine the first homology group of this Lagrangian is zero, this isotopy is exact.

Lemma \ref{rotate_lag} below tells us that fibers of $L_{\cup_j}$ and $L_{\cup_{j+1}}$ over each point of the diagonal intersect transversely at only one point.  We grade the $L_{\cup_{j}}$ in such a way that the absolute Maslov index of this intersection point (with regard to the two $S^2$ fibres) equals one.
  Construction for $\cap_{j}$ is similar.\\

In order for $\Phi$ to define a functor, we must verify that the above correspondences satisfy the same commutation relations as the tangles they are associated to. 
 The following is proven in \cite{RR1}
 
\begin{lemma}\label{rotate_lag}
If $|i-j|>1$ we have
\bqa L_{\cap_{i}}L_{\sigma_{j}}\simeq L_{\sigma_{j}}L_{\cap_{i}}  &\qquad& L_{\cup_{i}}L_{\sigma_{j}}\simeq L_{\sigma_{j}}L_{\cup_{i}} \\
\label{thesecondtolast}L_{\cap_{i}}\hspace{1ex} L_{\cup_{j}}&\simeq& L_{\cup_{j}}\hspace{1ex} L_{\cap_{i}} \quad \\
L_{\sigma_i}L_{\sigma_j}&\simeq&L_{\sigma_j}L_{\sigma_i}.\\
\eqa

 For any $i$ we have
 
 \begin{eqnarray}
L_{\sigma_i}L_{\sigma_{i+1}}L_{\sigma_i}&\simeq&L_{\sigma_{i+1}}L_{\sigma_i}L_{\sigma_{i+1}} \\
\label{2tolast}L_{\cap_{i}}L_{\sigma_{i}}\simeq L_{\cap_{i}}\{ 1\}   &\qquad& L_{\sigma_{i}}L_{\cup_{i}}\simeq L_{\cup_{i}}\{ 1\} \\
\label{last}L_{\sigma_{i\phantom{;}}} L_{\cup_{i+1}}\simeq L_{\sigma_{i+1}}^t L_{\cup_i }  &\qquad& L_{\sigma_{i+1}}L_{ \cup_{i}}\simeq L_{\sigma_{i}}^{t}L_{\cup_{i+1}}\\
\label{thelast}L_{\cap_{i;m}} \hspace{1ex} L_{\cup_{i+1;m}} &\simeq& L_{id_{m-1}} \{1\}.
\end{eqnarray}
 Here  ``$\simeq$'' means exact isotopy.

\end{lemma}

From \ref{yetterlemma} and the Lemma \ref{rotate_lag} we get the following.

\begin{theorem}\label{wellfunc}  The assignment $\Phi$  in \eqref{Phi} is a functor from the category of even tangles to the symplectic category. 
\end{theorem}

\subsection{ The symplectic Khovanov homology of even tangles}

 We can obtain a group valued  tangle invariant from the functor $\Phi$ as follows.

\begin{definition}
 \bqa\label{grpvalued}\khs(T)=\bigoplus_{\substack{C\in\pz{C}_m\\C'\in \pz{C}_n}}  HF(L_C^t,\Phi(T),L_{C'})\\
\label{grpvaluedCF} \cf(T)=\bigoplus_{\substack{C\in\pz{C}_m\\C'\in \pz{C}_n}}  CF(L_C^t,\Phi(T),L_{C'})
\eqa
\end{definition}
We will, in subsection \ref{vanishofdiff}, put extra conditions on the chain complex \eqref{grpvaluedCF} for $T$ a flat tangle.
Each summand in the above direct sum is equal to the Floer cohomology of the Lagrangians
\begin{eqnarray*}
\mathcal{L}_0&=&L_C\times L_0\times L_{1,2}\times...\times L_{2k-1,2k}  \\
\mathcal{L}_1&=&L_{0,1}\times L_{2,3} \times ...\times L_{2k+1}\times L_{C'}
\end{eqnarray*}
 in $\mathcal{Y}=\mathcal{Y}^-_{n_0}\times\mathcal{Y}_{n_1}\times ...\times \mathcal{Y}_n $. 
Since these Lagrangians are not necessarily compact one has to make extra effort to prove that the above Floer homology is well-defined.  This is done in \cite{RR1} by truncating the Lagrangians near infinity using the Stein structure on the manifolds $\mathcal{Y}_m$.

 \begin{theorem}\label{khswelldef}
 For any tangle $T$, $\khs(T)$ is well-defined and is independent of the decomposition of $T$ into elementary tangles.
  \end{theorem}

   It is clear that if $K$ is a $(0,0)$-tangle, i.e. a link, then the above invariant equals the original invariant of Seidel and Smith
  \eqref{defss}.


\section{Results on  symplectic Khovanov homology}\label{ResultsChapter}

\subsection{Maps induced by cobordisms}\label{minimalcobo}


In this subsection we study the maps induced on $\khs$ by  cobordisms.
%
%
%
Using elementary Morse theory one can decompose any cobordism $S$ between two  tangles $T,T'$ into elementary cobordisms
$$S=S_l\circ S_{l-1}\circ \cdots\circ S_1,$$ 
 where each $S_i$ belongs to one of the three elementary types discussed below.   
We assign homomorphisms to each one of these elementary types and so get  a homomorphism $\khs(S_i)$ for each $S_i$. 
To define these elementary cobordism maps one needs to decompose the tangles involved into elementary ones. However cobordism maps do not depend on such decomposition  of the tangle in an appropriate sense (Lemma \ref{cobocommdiag}).
 %
We then define $\khs(S)$ to be the composition 
$$\khs(S_l)\circ \khs(S_{l-1})\circ \cdots\circ \khs(S_1).$$ Since such a decomposition is not unique, one can potentially get different maps from different decompositions. We do \emph{not} address this problem here.
%
%


\textbf{\large{Type I.}}
Cobordisms, equivalent to trivial cobordism, between equivalent tangles. The (iso-)morphism assigned to such a cobordism is given by the functoriality theorem. C.f. Theorem \ref{khswelldef}.\\

\textbf{\large{Type II.}}
Birth or death of an unlinked circle:
\bqa
S_\bigcirc: \quad T\bigcirc T'\lra TT'\\
S_\bigcirc^t: \quad TT' \lra T\bigcirc T'.
\eqa
We know from \cite{RR1} that there is a canonical isomorphism 
$$\khs(T\bigcirc T')\cong \khs(TT')\otimes_{\Z} \pz{V}.$$
We define  the map 
\bq
\khs(S_\bigcirc): \khs(T\bigcirc T')\lra \khs(TT')
\eq
induced by the cobordism $S_\bigcirc$ to be $id\otimes \varepsilon$
 and the map 
\bq
\khs(S_\bigcirc^t):  \khs(TT')\lra \khs(T\bigcirc T')
\eq
to be $id\otimes \iota$. Here $\varepsilon$ and $\iota$ are the trace and unit maps from \ref{khovtqft}. Note that both maps can be defined as maps induced by quilts.

\begin{figure} \label{minimal_cobo_quilt}
\begin{center}
\includegraphics{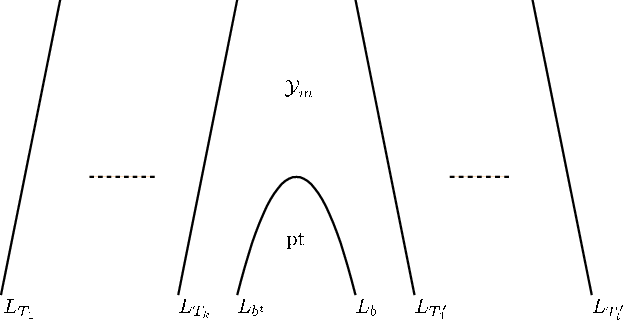}
 \caption{Quilt associated to a minimal cobordism}
\end{center}
\end{figure}

\textbf{\large{Type III.}}
Saddle point cobordisms:
\bqa
S_{\cap_i}:\quad T\cup_i\cap_i T' \lra TT'\\
S^t_{\cap_i}:\quad TT' \lra T\cup_i\cap_i T'.
\eqa
Let $\Phi(T), \Phi(T')$ be the generalized Lagrangian correspondences associated to $T$ and $T'$.
We  define the following cobordism maps as follows.
\bq\label{capcupcobo}
\cf(S_{\cap_i}):\cf(T\cup_i\cap_i T')\ra \cf(TT')
\eq
\bq\label{capcupcopro}
\cf(S_{\cap_i}^t):\cf(TT')\ra \cf(T\cup_i\cap_i T')
\eq
The homomorphism \eqref{capcupcobo} is defined to be the relative
invariant associated to the quilt in Figure \ref{minimal_cobo_quilt}. 
Recall \cite{WWquilts} that such a map is given by the count of the
zero dimensional part of the moduli space of \pse quilts as in Figure
\ref{minimal_cobo_quilt} with the indicated Lagrangian boundary conditions.
The homomorphism \eqref{capcupcopro} is the relative invariant of the transpose of this quilt.

\begin{figure} 
\begin{center}
\includegraphics{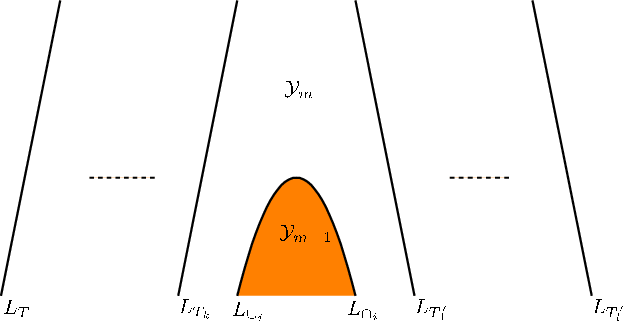}
 \caption{Quilt associated to a saddle cobordism}
 \label{minimalcobocap}
\end{center}
\end{figure}

Remember that we define the $\khs$ of a tangle by using a decomposition of it into elementary tangles and using different decompositions result in isomorphic abelian groups. The following lemma shows that elementary cobordism maps are natural with respect to change of decomposition (or, in other words, isotoping the  tangles involved).

\bl \label{cobocommdiag}
Let $T_0,T_1,T'_0,T'_1$ be tangles such that $T_i$ is equivalent to $T'_i$ for $i=0,1$. Then we have the following commutative diagram\\


\bd
\xymatrix{ \khs(T_0 \cup_i \cap_i T_1 ) \ar[d] \ar[r]^{\khs(S_b)}&       \khs(T_0\, T_1)\ar[d]\\
            \khs(T'_0 \cup_i \cap_i T'_1 )            \ar[r]^{\khs(S_b)}& \khs(T'_0 \,T'_1) }
\ed
where the vertical maps are isomorphisms. One has a similar diagram for the type II cobordisms.
\el
\bp
The vertical isomorphisms were constructed by showing that Lagrangian
correspondences assigned to elementary tangles satisfy the same
commutation relations as the corresponding tangles. So the maps are
given by Hamiltonian isotoping the corresponding Lagrangians and the
functoriality theorem. For the first kind, the Hamiltonian isotopy
induces a diffeomorphism between the corresponding moduli spaces of
quilts. The second kind is an instance of ``shrinking strips in
quilted surfaces '' and commutativity is given, in general settings,
by Theorem 5.1 in \cite{WWquilts}. The argument for type II cobordisms is similar.

\ep

\begin{figure} 
\begin{center}
\scalebox{1}{
\input{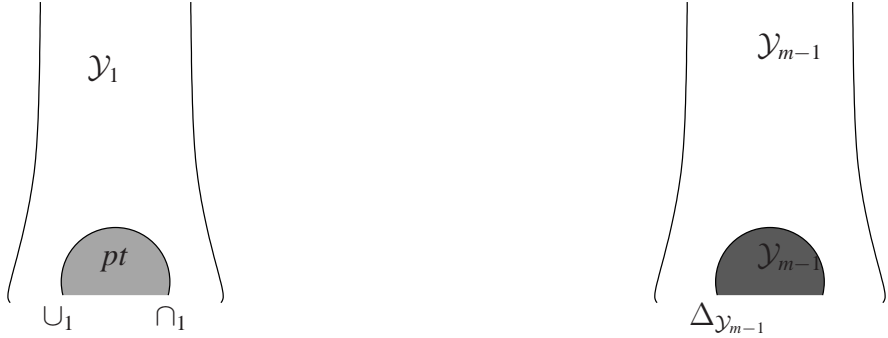}
}
\caption{Decomposing the quilt associated to a saddle cobordism}
\label{ConeQuiltDeCap}
\end{center}
\end{figure}

We next turn to the degree of saddle cobordism maps.
\bl
We have the following commutative diagram  where $HF(Q_{id})$ is the relative map of the quilt on the
right hand side of Figure \ref{ConeQuiltDeCap}.
\bd\label{degreediagram}
\xymatrix{ \khs(T\cup_{i}\cap_{i} T')\ar[d]^{\khs(S_{\cap_{i}})}  \ar[r]
  & \khs(T\, \on{id}\, \on{id}\, T')\otimes \khs(\cup_{1}
  \cap_{1})\ar[d]^{
HF(Q_{id}) \otimes\khs(S_{\cap_{1}})} \ar[r] & \khs(T T')\otimes \khs(\cup_{1} \cap_{1}) \ar[d]^{id\otimes\khs(S_{\cap_{1}})} \\
            \khs(T T')\ar[r] & \khs(T\, \on{id}\, \on{id}\, T' )\otimes \khs(\on{id}_1)
            \ar[r]& \khs(T T')\otimes \khs(\on{id}_1)
}
\ed
\el\label{degreelem}
\bp
We can use the local picture of the Lemma \ref{2nbhd}
to  isotope  the symplectic structure on $\Sm^{sub,\mu}\times
\mathbb{C}^3$ to the product one. This way  $L_{\cup_i}$ gets
smoothly isotoped to $ \Y_{m-1}\times S^2$. Floer homology is
invariant under such isotopies.  
Since
the Lagrangians are now Cartesian products and one can use product
complex structures to achieve transeversality, $\khs(S_{\cap_i})$ becomes
the tensor product of the relative maps of the quilts in Figure
\ref{ConeQuiltDeCap}. 
%
%
 The map $HF(Q_{id})$ is given by counting \pse strips of Maslov index zero and
so it is the identity map.
\ep

The degree of the vertical map on the right hand side is a priori zero (because it is a
pair of pants map )
 however because of the $-1$ degree shift coming from the extra cap in
 $T\cup_i\cap_i T'$ this degree equals one. 
Note that the horizontal maps on the right hand side are of degree 
$\frac{1}{2} \dim \mathcal{Y}_{m}$.
This is because the grading of $\Delta_{\mathcal{Y}_{m}}^t$ is by definition equal
to $\frac{1}{2}\dim \mathcal{Y}_{m}$ minus the grading on $\Delta_{\mathcal{Y}_{m}}$.
See also Corollary \ref{triangledecomp} below. 
In words the above commutative diagram implies that even though the
degree of a saddle cobordism map equals $\frac{1}{2} \dim
\mathcal{Y}_{m}$, when comparing to Khovanov homology this degree can
be taking to be one.


\textit{Minimal cobordisms} (c.f. Definition \ref{defminimalcobo}) are specific combinations of saddle cobordisms which deserve special attention.
 Let
$a\in \P_l, b\in\P_m$ and $c\in \P_n $.  
Recall from subsection \ref{ringshn} that we denote the minimal cobordism $b^tb\ra id$ by $S_b$. 
From $S_b$ we get the cobordism $\bf{1}_a\bf{1}_T  S_b \bf{1}_{T'}\bf{1}_c  $  between $aTbb^tT'c$ and $aTT'c$.
  We associated to this cobordism the quilt  depicted in Figure \ref{minimal_cobo_quilt}  and we denote it by $Q_b$. The relative invariant associated to $Q_b$ gives a homomorphism of chain complexes

  \bq\label{thehomoCF}   CF(Q_b): \quad CF(L_a^t,\Phi(T),L_b)\otimes_\Z CF(L_b^t,\Phi(T'),L_c)\longrightarrow CF(L_a^t,\Phi(T\circ T'),L_c). \eq

  as well as   a homomorphism of graded groups

\bq\label{thehomo}   \khs(S_b):=HF(Q_b): HF(L_a^t,\Phi(T),L_b)\otimes HF(L_b^t,\Phi(T'),L_c)\rightarrow HF(L_a^t,\Phi(T\circ T'),L_c). \eq
Summing over all such $a,b$ and $c$ we get maps
\bq
\cf(S_b): \quad \cf(T)\otimes_\Z \cf(T') \longrightarrow \cf(TT')
\eq
and
\bq\label{S_bmapH}
\khs(S_b): \quad \khs(T)\otimes_\Z \khs(T') \longrightarrow \khs(TT').
\eq
If $Q_b^t$ is obtained from $Q_b$ by reversing the incoming and outgoing ends we get

\bq\label{reversehomo}   \khs(S_b^t):  \khs(T\circ T') \rightarrow  \khs(\Phi(T))\otimes_\Z \khs(T'). \eq
It follows from the formula for the degree of the relative map of a quilt (Remark 2.4 in \cite{WWquilts}) that 
\bq
\deg \cf(S_b)=0= \deg \cf(S_b^t).
\eq

In the case that $T=T'=id_1$, $\mathcal{Y}_1$  is diffeomorphic to $T^*S^2$ and has a symplectic form which is cohomologous and homotopic to  the canonical symplectic structure on $T^{*}S^{2}$. So by Moser's trick, applied to a compact \nbhd of the zero section, it is  symplectomorphic to the standard one. Therefore  we get maps
\bq\label{mult}  m_{symp}=\khs(S_{\cap}): \khs(\bigcirc)\otimes\khs(\bigcirc)\to \khs(\bigcirc)\eq  
and
 \bq\label{comult} \Delta_{symp}=\khs(S^{t}_{\cap}): \khs(\bigcirc)\to \khs(\bigcirc)\otimes\khs(\bigcirc).\eq


\bl\label{m=msymp}

There is an isomorphism $\phi: \khs(\bigcirc)\to \pz{V}$ under which $\Delta_{symp} $ and $m_{symp}$ correspond to $\Delta$ and $m$ respectively.

\el
\bp
Let $f_0,f_1,f_2$ be three Morse functions on a Riemannian manifold $M$.
Fukaya and Oh \cite{FukayaOh} 
 prove that if  we equip  the  cotangent bundle of $M$ with the almost complex structure induced by the Levi-Civita connection on $M$ then, for generic choice of the $f_i$,
there is an orientation preserving diffeomorphism between the
moduli space of pseudoholomorphic triangles   connecting intersection points of the $df_i$ and the moduli space of pair of pants trajectories between the corresponding critical points of $F=f_0-f_1, G=f_1-f_2, H=f_2-f_0$.
The moduli space of \pse triangles in $T^*S^2$ with the almost complex structure from $\mathcal{Y}_1$ is zero dimensional and cobordant to the moduli corresponding to the almost complex structure induced by Levi-Civita connection. So the sum of the elements of the two are equal. 
Therefore after choosing (unique) homogeneous generators $1$ and $X$ for $\khs(\bigcirc)\cong H^{*}(S^{2})$, \eqref{mult} corresponds to the wedge product on $\pz{V}$ which is in turn equal to $m$. 

The same arguments  show that $\Delta_{symp}$ corresponds to the operation given by counting inverted Y's in Morse homology of $S^{2}$.
A direct computation shows that this operation is equal to $\Delta$.

\ep

\begin{figure}[ht]\label{quilt_associativity}
\begin{centering}
 \includegraphics[scale=.6]{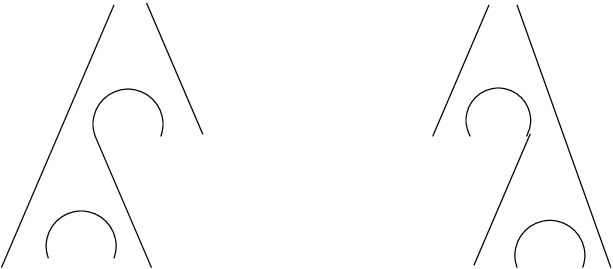}
 \caption{Isotopy between the composition of two quilts}
\end{centering}
 \end{figure}

Note that the horizontal composition $SS'$ of two cobordisms $S$ and $S'$ equals $(S\; id)\circ (id\;  S')$.
\bl
For two minimal cobordisms $S_a$ and $S_b$ we have
 \bq\label{horizcommut} \khs(S_a\; id)\circ \khs(id\; S_b)=\khs(id\; S_b)\circ \khs(S_a\; id).\eq
\el
\bp
From \cite{WW} we know if $Q$ and $Q'$ are two quilts which can be composed vertically (i.e. along the strip-like ends) then  $HF(Q\circ Q')=HF(Q)\circ HF(Q')$. We also know that $HF(Q)$ is (at the cohomology level) invariant under the isotopy of the quilt $Q$. The lemma follows from these two facts together with the isotopy in Figure \ref{quilt_associativity}.
\ep

\subsection{$H^m$ module structure}\label{Hnmodulestr}

We define the symplectic analogue of the rings $H^n$ as \bq\label{Hsymp}H^n_{symp}=\khs(id_n)=\bigoplus_{a,b\in \P_n} HF(L_a^t,L_b).\eq 
 The product map from $HF(L_a^t,L_b)\otimes HF(L_c^t,L_d)$ to $HF(L_a^t,L_d)$ is  zero if $b\neq d$ and is given by the map  $HF(Q_b)$ otherwise.

\bl\label{HFcrossingless} For any $a,b\in \P_m$ we have $HF(L_a^t,L_b)\cong H^*(S^2)^{\otimes k}\cong \pz{V}^{\otimes k}\{k\}$ where $k$ is the number of circles in $a^t b$.
\el
\bp  The Lagrangian $L_a$ equals the composition of the Lagrangians associated to its arcs and similarly for $L_b$ so $\Phi(a^t)\#\Phi(b)=\Phi(a^tb)$. Therefore
\bqas
 HF(L_a^t,L_b)&=&HF(\Phi(a^t),\Phi(b))\{m\}=HF(\Phi(a^tb))\{m\}\\
&=&\khs(\Phi(k\bigcirc))\{m\}=H^*(S^2)^{\otimes k}. 
\eqas
\ep

\begin{theorem}\label{HsympH} For each $n$ we have an isomorphism $$H^n_{s}\cong H^n$$ as graded algebras over $\Z/2$.
\end{theorem}
\bp

Let $a,b,c \in \P_m$. Observe that\\ $HF(\Phi(a^t),\Phi(b))\otimes HF(\Phi(b^t),\Phi(c))$ is canonically isomorphic to $HF(\Phi(a^t\, b\, b^t\, c))$. Label the arcs of $b$  with numbers $1$ to $m$. The Lagrangian correspondence $\Phi(a^t \, b\, b^t\, c)$ is equivalent in the symplectic category to $\Phi(C_1)$ where $C_1$ is unlinked disjoint union of some $k$ circles. Each of these circles has some marked points on it corresponding to the arcs of $b$ and $b^t$ and these marked points are labelled  with numbers which specify the cobordisms in $S_b: a^t\, b\, b^t\,  c\ra a^t\, c $.
 So  $S_b$ is equivalent to the composition of cobordisms $S_i: C_i\ra C_{i+1}$, $i=1,\ldots m$, where  each $C_i$  is an unlinked disjoint union of a number of circles and each $S_{i}$ is either a saddle cobordism or it permutes two circles. 
We observe that each $\Phi(C_i)$ is a generalized Lagrangian correspondence  of the form
\bq\label{Hmlag}
pt \ra \mathcal{Y}_1 \ra pt \ra \mathcal{Y}_1 \ra pt \ra\ldots \ra pt \ra \mathcal{Y}_1 \ra pt.
\eq
 Therefore we can use the Lemma \ref{m=msymp} to conclude that we have a commutative diagram
\bd
\xymatrix{
\khs(C_i)\ar[r]^{\khs(S_i)} \ar[d] &\khs(C_{i+1})\ar[d]\\
\kh(C_i)\ar[r]^{\kh(S_i)} & \kh(C_{i+1})
}
\ed
where (because of the special form of the Lagrangians \eqref{Hmlag}) the vertical arrows are canonical isomorphisms.


Now Lemma  \ref{cobocommdiag} ensures that the isomorphisms $\khs(a^t\, b\, b^t\, c)\cong \khs(C_1)$ and $\khs(a^t c)\cong \khs(C_m)$ intertwine the cobordism maps 
$$\khs(S_b):\khs(a^t\, b\, b^t\, c)\ra \khs(a^t c)$$ and $$\khs(S_m)\circ \cdots\circ \khs(S_1): \khs(C_1)\ra \khs(C_m).$$


\ep
For an $(l,m)$-tangle $T$, $\khs(T)$ has a structure of a $(H_{s}^l,H_{s}^m)$ bimodule as follows.
We have $$\khs(T)=\bigoplus_{b\in\P_l,c\in \P_m} HF(L_b,\Phi(T),L_c).$$
The part $_a H^l_b$ of $H^{l}_{s}$ acts on $HF(b,\Phi(T),c)$ from left by the map $\khs(S_b)$ (in \eqref{S_bmapH}). So does $_c H_d^m$ from right by the map $\khs(S_c)$. We set the left action of $_a H^l_{b'}$ on $HF(b,\Phi(T),c)$ to be zero if $b\neq b'$ and similarly for the right action.
This way we obtain an $(H^l,H^m)$-bimodule structure on $\khs(T)$. 
 
\begin{remark}\label{cssismodule}
Note that since the cobordism maps $\cf(S_b)$ are of degree zero, the chain complex $\cf(T)$ can be regarded as a chain complex of $(H^l,H^m)$-bimodules.
\end{remark}

\bl\label{eqforelementangle}
With the same notation as in \eqref{HmHm-1} we have
\bqas
\khs(id_m)&=&H^m=\pz{Kh}(id_m)\\
\khs(\cap_{i;m})&=&\pz{Kh}(\cap_{i;m})\\
\khs(\cup_{i;m})&=&\pz{Kh}(\cup_{i;m})\\
\label{khsofbraid}
\khs(\sig_{i;m}^\pm)&=&\overline{\kh}(\sig_{i;m}^\pm)= \\
&&\left(H^{m-1}   \oplus {H^m}' \oplus {H^m}'\oplus H^m_1 \right)\otimes\pz{V}\{1\mp2\} \bigoplus H^m_2\otimes \pz{V}\{2\}
\eqas
as $H^m$ modules.
\el
\bp
The first three equations follow from the fact that  $\kh$ and $\khs$  for disjoint union of $k$ circles are equal to $\pz{V}^{\otimes k}$. For the las one, the first equality was proved in \cite{RR1} 
 and the second equality in the lemma \ref{khofbraid}.

\ep


\subsection{Vanishing of the differential for flat tangles}\label{vanishofdiff}

\bl\label{crossmdiff0}
Let  $C_1,C_2\in\P_m$ be two crossingless matchings. Then we can choose Floer data in such a way that the Floer chain complex $CF(L_{C_1},L_{C_2})$ has differential equal to zero.
\el
\bp
We prove by induction on $m$. If $m=1$ then there is only one crossingless matching and the Floer chain complex equals $CF(S^2,S^2)$ where $S^2$ is the zero section in $\mathcal{Y}_m=T^* S^2$. We can Hamiltonian isotope  the zero section to a Lagrangian $L$ s.t. $L$ intersects the zero section at only two points. For example we can take $L$ to be the graph of the one-form $df$ where $f$ is the height function on the zero section.  In this case the Floer differential has to be zero because otherwise $HF(S^2,L)$ will not be equal to $H^*(S^2)$. This can also be seen by considering the Maslov indices of intersection points.

Now assume the statement holds for all crossingless matchings in $\P_k$ for $k< m$. Let $\alpha_1$ be an arbitrary arc in $C_1$ and $\mu_1,\mu_2$ its endpoints. There are two cases. Either there is an arc $\alpha_2$ in $C_2$ joining $p$ and $q$ or there is no such arc.  Proof for these two cases are similar to the proofs of the Kunneth formula and the Thom isomorphism for Floer homology \cite{SS}.
In the first  case let $\bar{C}_i$  be obtained from $C_i$ by deleting   $\alpha_i$, $i=1,2$. Then we can use lemma \ref{2nbhd} and then isotope the induced metric into the product metric. So $L_{C_i}$ gets isotoped to $L_{\bar{C}_i}\times S^2$. We choose a time dependent almost complex structure $\bar{J}_t$ on the base which is a compactly supported perturbation of its standard structure $\bar{J}_0$. We choose the almost complex structure on the total space to be equal to the product $\bar{J}_t\oplus J_{\mathbb{C}^3}$ in a small \nbhd  $U_0$ of the zero section and equal to $\bar{J}_0\oplus J_{\mathbb{C}^3}$ outside an open set $U_1$ containing $\bar{U_0}$. This way we can obtain an almost complex structure which is both regular and has similar properties to the product structure inside $U_0$. Since  our pseudoholomorphic strips are confined to $U_0$,
 we have
 $$CF(L_{C_1},L_{C_2})= CF(L_{\bar{C}_1},L_{\bar{C}_2})\otimes CF(S^2,S^2). $$
So the claim follows from the induction hypothesis and the argument for the base case.  In the second case let $\alpha_2 $ be the unique arc in $C_2$ which has $\mu_2$ as an endpoint and let $\mu_3$ be its other end point.  Now we can use lemma \ref{3nbhd} to identify $L_{C_i}$ with $L_{\bar{C}_i}\times_{S^1} \Lambda_{\alpha_i}$ where $\Lambda_\alpha$ is the Lagrangian sphere associated to the curve $\alpha$ as defined in subsection \ref{fibered}. We choose the almost complex structure in a  way similar to that of the first case above.  Note that there are two possible configurations of the curves $\alpha_i$. %

In either case $\Lambda_{\alpha_1}$ and $\Lambda_{\alpha_2}$ intersect at only one point  $p$ corresponding to $\mu_2$. So we have $$CF(L_{C_1},L_{C_2})=CF(L_{\bar{C}_1},L_{\bar{C}_2})\otimes \bb{Z}<p>.$$
Let $u$ be a holomorphic strip joining to intersection points of $L_{C_1}$ and $L_{C_2}$. So we have $u=(u',u'')$ where $u'$ is the projection to the first factor. By the induction hypothesis, $u'$ is constant.  Projection to the second factor is a holomorphic strip in $\bb{C}$ which has its boundary on $\alpha_1$ and $\alpha_2$. Such a finite energy curve has to be constant  by the exponential convergence property of pseudoholomorphic strips.  Therefore $u''$ is also constant so we get the desired result.

\ep

\bl\label{flatdiff0}
Let $T$ be a flat $(m,n)$-tangle. 
We can choose the Floer data in such a way that the Floer chain complex whose cohomology is $\khs(T)$ has differential equal to zero.
\el
\bp
Let $T=T_1\cdots T_{k-1}$ be a decomposition of $T$ and let $T_0\in \P_m$ and $T_k\in \P_n$. Let $L_{T_i}$  is a correspondence between $\mathcal{Y}_{m_i}$ and $\mathcal{Y}_{m_{i+1}}$. We use induction on $m=\sum m_i$. The case $m=1$ was treated in Lemma \ref{crossmdiff0}. If $T_1$ is the identity tangle then $CF(L_{T_0}, L_{T_1},\cdots, L_{T_k})=CF(L_{T_0}, L_{T_2},\cdots L_{T_k})$. So we can assume that $T_1$ is a cup. Therefore the both $L_{T_0}$ and $L_{T_1}$ are obtained by relative vanishing cycle construction from Lagrangians in $\mathcal{Y}_{m_0-1}$ and $\mathcal{Y}_{m_0-1}\times\mathcal{Y}_{m_1}$. Therefore we can use the same argument as in the proof of \ref{crossmdiff0} for the induction step.
\ep

\bdf\label{CFflat} For a flat $(m,n)$-tangle $T$ we require the chain complex 
$$\cf(T)=\bigoplus_{a\in \P_m, b\in \P_n} CF(L_a^t,\Phi(T),L_b)$$
to be given by Floer data in lemma \ref{flatdiff0}. 
\edf


\subsection{Composition property of $\khs$ for flat tangles}\label{functorialityforflat}
Let $T$ and $ T'$ be  $(l,m)$ and $(m,n)$ tangles respectively.
Consider the map $\psi_s$

\bd\label{tensormap}
\xymatrix{\khs(T)\otimes_{\bb{Z}}\khs(T')\ar@{=}[d]\\
  \bigoplus_{a,b,b',c} HF(L_a^t,\Phi(T),L_b) \otimes_{\bb{Z}}HF(L_{b'}^t,\Phi(T'),L_c)\ar[d]^{\psi_s}\\ 
    \bigoplus_{a,c} HF(L_a^t,\Phi(T'\circ T), L_c)\cong \khs(T'\circ T)
}
\ed
\linebreak
which is zero if $b\neq b'$ and equals $\khs(\mbf{1}_a \mbf{1}_T S_b  \mbf{1}_{T'} \mbf{1}_c)$ otherwise. Here, as before,  $S_b$ is the minimal cobordism between $bb^t$ and $id_{m}$.
The abelian group $\khs(T)\otimes_{\bb{Z}}\khs(T')$ has the structure of a $(H^l,H^n)$-bimodule and $\psi$ is a $(H^l,H^n)$-bimodule map.
If $x\in HF(L_a^t,\Phi(T),L_b)$, $y\in  HF(L_b'^t,\Phi(T'),L_c)$ and $\xi\in\phantom{j}_{b}H^m_{b'}$ then

\bqas
\psi_s(x\xi \otimes y )&=&\khs(\mbf{1}_a \mbf{1}_T    S_{{b'}^t} \mbf{1}_{T'} \mbf{1}_c)\khs(\bf{1}_a\bf{1}_TS_b\bf{1}_{b'}\bf{1}_{T'}\bf{1}_c) (x,y)\\
\psi_s(x\otimes \xi y)&=&\khs(\mbf{1}_a \mbf{1}_T    S_{{b}} \mbf{1}_{T'} \mbf{1}_c) \khs(\bf{1}_a\bf{1}_T \bf{1}_{b} S_{{b'}^t} \bf{1}_{T'}\bf{1}_c)(x,y)    \eqas
It follows from \eqref{horizcommut} that these two are equal and so $\psi_s $ factors through a map of bimodules $\khs(T)\otimes_{H^m}\khs(T')\rightarrow \khs(T\circ T') $ which we still denote by $\psi_s$.

\bpr\label{tensorkhsymp}
If $T$ and $T'$ are flat then
$\psi_s$ gives an isomorphism

 $$\psi_s:\khs(T)\otimes_{H^m}\khs(T') \cong\khs(T\circ T').$$

\epr

\bp Proof is exactly the same as that of Theorem 1 in \cite{functorvalued}.
 The map $\psi $ is the direct sum of the maps
 $$_a\psi_c:\oplus_b HF(L_a^t,\Phi(T),L_b)\bigotimes_{H^m}  \oplus_{b'}HF(L_{b'}^t,\Phi(T'),L_c)\rightarrow HF(L_a^t,\Phi(T\circ T'),L_c)$$

We have $\oplus_b HF(L_a^t\Phi(T),L_b) \cong \khs(a^t T)$ and 
$$\bigoplus_{b'}HF(L_{b'}^t,\Phi(T'),L_c)\cong \khs(T'c )\{n\}$$
 as left and right $H^m$-modules respectively. We also have $HF(L_a^t,\Phi(T\circ T'),L_c)\cong \khs(aTT'c)\{n\}$. Therefore the argument is reduced to showing that
$$ \khs(a^tT)\otimes_{H^m} \khs(T'c)\cong \khs(a^tTT'c).$$

Now $a^tT$ and $T'c$ are $(0,m)$ and $(m,0)$-tangles respectively so $a^tT=a'\oplus i\bigcirc$ and $T'c=c'\oplus j\bigcirc$ where $a',b'\in\P_m$ and
$i$ and $j$ are the number of circles in $a^tT$ and $T'c$ respectively. Thus we have $\khs(a^tT)=\khs(a')\otimes \pz{V}^i$, $\khs(T'c)=\khs(c')\otimes \pz{V}^j$ and $\khs(a^t TT'c)=\khs(a'c')\otimes \pz{V}^{i+j}$. So we need to show that $$\khs(a')\otimes_{H^m} \khs(c')=\khs(a'c').$$

We have $H^m\otimes_{H^m} H^m=H^m $ and if we multiply this identity with the idempotent $1_{a'}$ from left and by $1_{c'}$ from right we get the desired result. 
\ep

\begin{corollary}\label{khs=khflat}
For any flat tangle $T$ we have $$\khs(T)=\overline{\kh}(T)=\kh(T).$$
\end{corollary}
\bp This follows from \eqref{khdef}, \ref{tensorkhsymp} and \ref{eqforelementangle}.\ep

\bpr\label{coboequflat}
Let $T,T'$ be flat $(m,n)$-tangles and $S$ a cobordism between $T$ and $T'$ which equals a composition of minimal cobordisms.   Then, with coefficients in $\Z/2$, we have
\bq
\khs(S)=\kh(S).
\eq
\epr
\bp
By \eqref{horizcommut} we can assume that $S$ consists of  a single minimal cobordism.  Therefore we have
$T=T_1 c c^t T_2$  and $T'=T_1 id \, T_2$  for a crossingless matching $c$ and $S$ equals $\mbf{1}_{T_1}S_c\mbf{1}_{T_2}$. For any $a\in \P_m$ and $ b\in \P_n$,
$a^t T_1$ equals  a crossingless matching   $a_1^t \in \P_m$ disjoint union with some $k$ circles. The same is true for
$T_2 b$ i.e. $T_2 b$ equals $b_2\in \P_n$ disjoint union with $l$ circles.
So, the problem is  reduced to showing that the map
 $$\khs(S)=HF(Q_c):HF(L_{a_1}^t,L_c,L_c^t,L_{b_2})\otimes \pz{V}^{k+l} \longrightarrow HF(L_{a_1}^t,L_{b_2})\otimes \pz{V}^{k+l}$$
equals $\kh(S)$. But 
$$HF(L_{a_1}^t,L_c,L_c^t,L_{b_2})=HF(L_{a_1}^t,L_c)\otimes HF(L_c^t,L_{b_2})=\, _{a_1}H_c\,\otimes \,_cH_{b_2}$$ and $HF(L_{a_1}^t,L_{b_2})= \,_{a_1}H_{b_2}.$ Therefore the lemma follows from the isomorphism  of the algebra structures on $H^m$ and $H^m_{s}$ over $\Z/2$ (Lemma \ref{HsympH}).


\ep

\subsection{ Skein exact triangle for symplectic Khovanov homology}\label{exacttri}
In this subsection we prove an exact triangle for the Seidel-Smith invariant which is analogous to skein relations for knot polynomials.  
%
%
%
The tool we use is the exact triangle for Lagrangian Floer homology. 
This exact triangle was discovered by Seidel \cite{seideltriangle} for Dehn twists. We use a generalization of this triangle to fibred Dehn twists due to Wehrheim and Woodward \cite{WWtriangle}.
%
Let $M$ be a symplectic manifold and  $C\subset M$ a spherically fibred coisotropic fibring over a base $B$. We denote the fibred Dehn twist along $C$ by $\tau_C$. The embedding $(\iota\times\pi)C$ is a Lagrangian submanifold of $M^-\times B$.
By the abuse of notation we sometimes denote this submanifold by $C$.

\begin{figure}[ht]
\begin{center}
\scalebox{0.6}{
\input{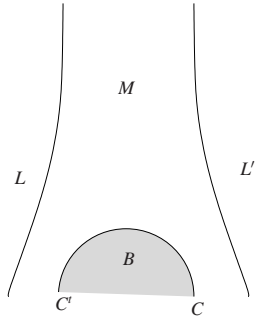}
}
\caption{The quilt used in the exact triangle}
 \label{conequilt}
\end{center}
\end{figure}

Let $Q_0$ be the quilt in the Figure \ref{conequilt}. The exact triangle in \cite{WWtriangle} establishes a  quasi-isomorphism between the Floer chain complex $CF( L,\tau_C L')$ and the cone of the morphism
$f:=CF(Q_0)$, i.e. the relative map associated with $Q_0$.

\bq\label{conesmap}
f=CF(Q_0):\quad CF(L,\; (\pi\times\iota)C^t, (\iota\times\pi)C,\; L')\{-\half\dim B\} \longrightarrow CF(L,L').
\eq

More precisely we have the following.

\begin{theorem}[Wehrheim, Woodward \cite{WWtriangle}]\label{trithmWW}
If $C\subset M$ has codimension at least two and the triple ($L_0,L_1,C)$ is monotone and has Maslov index greater than or equal 3 then there is a quasi-isomorphism $(h\{1\},k)$ from $$Cone(f)=CF(L,C^t,C,L')\{-\half\dim B+1\}\bigoplus CF(L,L')$$ to
$CF( L,\tau_C L')$.
\end{theorem}



Recall that a Lagrangian correspondence $L\subset M^{-}\times N$ between two symplectic manifolds  is proper if for each $y\in N$ the set $\{x\in M| (x,y)\in L\}$ is compact. (For Lagrangian submanifolds this is equivalent to compactness.) We call  a Lagrangian submanifold $L$ of a Stein manifold $(M,\psi)$ \emph{allowable} if the critical point set of $\psi|_{L}$ is compact. When doing Floer theory in a Stein manifold we use compatible almost complex structures which are invariant under the Liouville flow outside a compact subset.

\bpr\label{trithmCF}
If $(M,\psi)$ is Stein,  $C\subset M$ has codimension at least two and  each one of $L_0,L_1,C$  is exact, proper and allowable  
then the conclusion of Theorem \ref{trithmWW} holds. 
\epr

\bp
Properness implies that the intersection in $(L,C^{t},C,L')$
(i.e. $L\times C \cap C^{t}\times L'$) is compact. In \cite[Lemma
3.3.2]{RR1} it was shown that for any two allowable Lagrangian
submanifolds $\L_{0},\L_{1}$ of a Stein manifold $(M,\psi)$, any \pse
curve with boundary on $\L_{0}\cup \L_{1}$ lies in $M_{\psi\leq C}$
where $C$ is the maximum of $\psi$ on $\L_{0}\cap \L_{1}$. This
implies that the moduli spaces of \pse quilts used in the statement
and  proof of Theorem \ref{trithmWW} are compact in our case.  The high Maslov index assumption in \ref{trithmWW} is to rule out bubbling; in our case bubbling is ruled out by exactness of the Lagrangians.    

\ep

At the  \aifty  level one has the following exact triangle in $\DGF(M)$.  
\begin{displaymath}\label{thecone}
\xymatrix{
& \gra \tau_{C}  \ar[dl]  &\\   C^t\#C\{-\frac{1}{2} \dim B\} \ar[rr] &&  \Delta_M\ar[ul]
 }
\end{displaymath}
Here  $\DGF(M)$ is the generalized Fukaya category of a Stein manifold as described in \cite[Section 4.3]{RR1}.
The objects of $\DGF(M)$ are proper exact allowable generalized Lagrangian submanifolds of $M$. 
This category is somewhat similar to wrapped Fukaya category \cite{AbSeidel}. The difference is that in the wrapped Fukaya category one uses the Reeb Flow instead of the Liouville flow and also one takes the direct limit of the Floer chain complexes of the images of the Lagrangian submanifolds under this flow. In our case the more restrictive properness assumption frees us from taking direct sums. 

If $\ul{L}=(L_k,L_{k-1},\ldots,L_1)$ is any generalized Lagrangian submanifold of $M$ then by applying the $A_\infty$ functor $\Phi_{\ul{L}}^\#=\Phi_{L_k}^\#\circ \cdots \circ\Phi_{L_1}^\#$ to \eqref{thecone} we get the following exact triangle in $\DGF(M)$.
\begin{displaymath}
\xymatrix{
& \gra \tau_{C} \# \ul{L} \ar[dl]  &\\   C^t\#C\#\ul{L}\{-1/2 \dim B\} \ar[rr] &&  \ul{L}\ar[ul]
 }
\end{displaymath}
%
%
Therefore theorem \ref{trithmCF} holds, without any change, if $L,L'$ are generalized Lagrangian submanifolds of $M$. One can prove this fact without using Fukaya categories. 



With the same assumptions as in \ref{trithmCF} let $M=M_1\times M_2$, $B=B_1\times M_2$,  and
 $C$ be of the form $C_1\times M_2$ where $C_1$ is a sphere bundle over $B_1$. Further assume that there are Lagrangian submanifolds $L_i\subset M_i$ and $L'_i\subset M'_i$ for $i=1,2$ \st
$L=L_1\times L_2$ and $L'=L'_1\times L'_2$. Let $d=-\frac{1}{2} \dim B$  and $d_1=-\frac{1}{2} \dim B_1.$
%
%
 Consider the map

\bq
CF(\bar{Q}_0):\quad CF(L_1,\; C_1^t, C_1,\; L'_1)\{-\frac{1}{2}\dim B_1\} \longrightarrow CF(L_1,L'_1).
\eq
\bcor\label{triangledecomp}

 $CF( L,\tau_C L')$ is quasi-isomorphic to 
\bq
Cone(\bar{Q})\otimes CF(L_2,L'_2)
\eq
and we have a commutative diagram
\begin{displaymath}
\xymatrix{
 CF(L,\; C^t, C,\; L')\{d\}    \ar[r]^{CF(Q_0)}\ar[d]^{\{\frac{1}{2}\dim M_2\}} &  CF(L,L')\phantom{CCC} \ar[d]\\
 CF(L_1,C_1^t, C_1, L'_1)\{d_1\} \otimes CF(L_2,L'_2)    \ar[r]^{CF(\bar{Q}_0)} &   CF(L_1,L'_1)\otimes CF(L_2,L'_2)\phantom{CCCC,C}
  }
\end{displaymath}

\begin{figure}[ht]
\begin{center}
\scalebox{1}{
\input{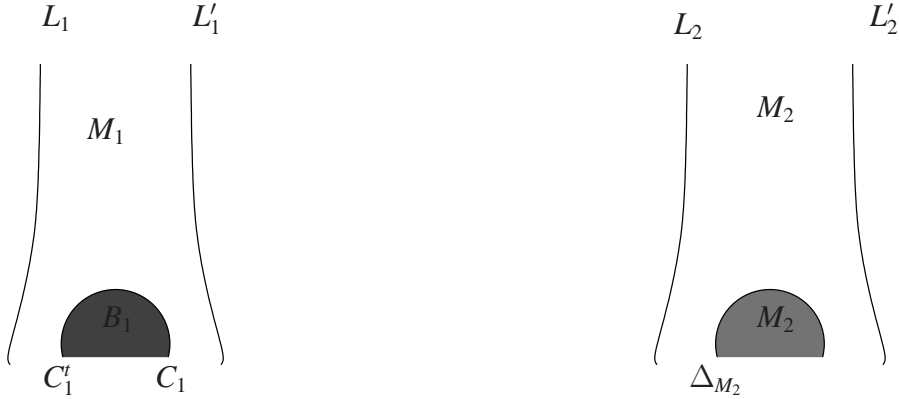}
}
\caption{Decomposition of the cone in Corollary \ref{triangledecomp}}
 \label{ConeQuiltDecomp}
\end{center}
\end{figure}

\ecor
\bp
We observe that $\tau_C=\tau_{C_1}\times id_{M_2}$. One can choose
product complex structures on $M_{1}\times M_{2}$ to achieve
transeversality.  So the homomorphism $CF(Q)$ is isomorphic to the tensor product of the maps induced by the two quilts in the Figure \ref{ConeQuiltDecomp}. The quilt on the left is $\bar{Q}$ and the quilt on the right induces the identity map. 
There is a grading shift
$CF(L_2,\Delta_{M_2}^t,\Delta_{M_2},L'_2)=CF(L_2,L'_2) \{\frac{1}{2}\dim M_2\}$
  coming from the fact that the grading of $\Delta_{M_2}^t$ equals $\frac{1}{2}\dim M_2$ minus the grading on $\Delta_{M_2}$.
\ep

\bpr\label{triforgeneralized}
With the same assumptions as in Theorem \ref{trithmCF}, let $\ul{L},\ul{L}'$ be two generalized Lagrangian submanifolds  of $M$. Then $CF(\ul{L}^t, \on{graph(\tau_C)}, \ul{L}')$ is quasi-isomorphic to the cone of the map
\bq
CF(\ul{L}^t, C^t,C, \ul{L}')\{-1/2\dim B   \}   \longrightarrow CF(\ul{L}^t, \ul{L}')
\eq 
\epr
\bp
Let $\ul{L}=(L_n,\cdots,L_k)$ and $\ul{L}'=(L_{k-1},\cdots, L_1)$ where $L_i\subset M_{i+1}^-\times M_{i} $ and $M_k=M$.  We can assume, by adding identity Lagrangian correspondences if necessary, that 
$$\L_0=L_n\times L_{n-2}\times\cdots \times L_k\times L_{k-1}\times\cdots\times L_2$$ and
$$\L_1=L_{n-1}\times L_{n-3}\times\cdots \times \on{graph}(\tau_C)\times\cdots\times L_1.$$

We have $\L_1=\tau_{C'}\left( L_{n-1}\times L_{n-3}\times\cdots \times L_1  \right)$
where $C'= M_{n+1}\times \cdots  \times M_{k+1}\times C\times M_{k-1}\times\cdots \times M_1  $  which fibers over $M_{n-1}\times\cdots   \times B\times \cdots\times M_{1}$.
The result follows from \ref{triangledecomp} by taking $M_1= M=M_k$ and $M_2$ to be the product of the rest of the manifolds $M_i$.

\ep

\begin{corollary}\label{triadjoint}
 $CF(L,\tau_C^{-1} L')$ 
is quasi-isomorphic to the cone  $\on{Cone}\left(CF(Q_0^t)\right)\{-1\}$.
\end{corollary}

\bp This is a standard argument.  
If  $l=\dim L$  then  we have \\
$CF^*(L,\tau_C^{-1} L')=   CF^*(\tau_C L,L')=  CF^{l-*}( L', \tau_C L)^\vee=\Hom(CF^{l-*}(L',\tau_C L),\mathbb{Z}  ).$
It follows from \ref{trithmCF}  that this is quasi-isomorphic to 
%
$$CF^{l-*}(L',\; (\iota\times\pi)C, (\pi\times\iota)C^t,\; L)\{1\}^\vee\oplus CF^{l-*}(L',L)^\vee$$
(with appropriate differential).
This  in turn equals
$$ CF(L,L')\bigoplus CF(L,\; (\iota\times\pi)C^t, (\pi\times\iota)C^t,\; L')\{-1\}.$$

\ep


Now we use Theorem \ref{trithmCF} (or more precisely \ref{triforgeneralized}) to obtain an exact triangle for the Seidel-Smith invariant.
In the case under study  $M=\mathcal{Y}_{l}$ and the spherically fibred isotropic is  $C=L_{\cap_i}$. The coisotropic submanifold $ C$ is a sphere bundle over  $B=\mathcal{Y}_{l-1}$. 
The fibred Dehn twist $\tau_C$ along $C$ equals the monodromy map $h_{\sig_i}$ and so for any Lagrangian $L\subset M$ we have
\bq
 L_{\sig_i} \circ L\simeq\tau_C L.
\eq
%
Let $kTl$ and $lT'm$ be tangles, $\sig^+_i, \sig^-_i\in Br_{2l} $  elementary braids and $T^\pm=T\sig^\pm_i T'$.
 We observe that if $a$ and $b$ are crossingless matchings and we take $\ul{L}=(a^t,\Phi(T))$ and $   \ul{L}'=(\Phi(T'),b^t)$ then the map $f$ in \eqref{conesmap} is the same as the cobordism map \eqref{capcupcobo}.

\bthm\label{SSexacttri}

Let $e$ be the difference between the number of negative crossings in $T^\pm$  
and $T\cup_i\cap_i T'$ (with the latter oriented arbitrarily).
 Then
$\cf (T^-)$ is quasi-isomorphic to the cone of

\bq\label{SSptri}
\cf (T\cup_i\cap_i T')\{-2e\} \stackrel{CSS(\mbf{1}_T S_{\cap_i}\mbf{1}_{T'})}\longrightarrow \cf (T T') 
\eq
 and $\cf(T^+) $ is quasi-isomorphic to the cone of

\bq\label{SSntri}
\cf (T T')\{-1\} \stackrel{CSS(\mbf{1}_T S^t_{\cap_i}\mbf{1}_{T'})}\longrightarrow  \cf (T\cup_i\cap_i T')\{  -1-2e\}.
\eq

\ethm

\bp
We note that the sign conventions for positive braids and positive  Dehn twists are opposites of each other.
%
Since the degree of the map \eqref{capcupcobo} equals $1$, we apply the degree shift $\{-1\}$ to its target to obtain a map of degree zero.  We have 
$w(T^-)=w(T\cup_i\cap_i T')+2e$  and $w(T^-)=w(TT')-1$
so we obtain \eqref{SSptri}.


In this case of $T^+$ we  have  $w(T^+) =w(T\cup_i\cap_i T')+2e$, $w(T^+)=w(TT')+1$ and   the cobordism map is of degree $1$.
\ep

\bcor
We have the following exact triangles

\begin{displaymath}\label{thecone}
\xymatrix{
&  \cf(T^-)  \ar[dl]_{\{1\}}  &\\   \cf(T\cup_i\cap_i T')\{-2e\} \ar[rr] &&  \cf(TT')\ar[ul]
 }
\end{displaymath}

\begin{displaymath}\label{thecone}
\xymatrix{
&  \cf(T^+)  \ar[dl]_{\{1\}}  &\\   \cf(TT')\{-1\} \ar[rr] &&  \cf(T\cup_i\cap_i T')\{-1-2e\}\ar[ul]
 }
\end{displaymath}

which give the following cohomology exact sequences.

\bq\label{khsm}
\cdots \ra \khs^{i-1}(T^-)\ra \khs^{i-2e}(T\cup_i\cap_i T') \ra \khs^i(TT')\ra \khs^{i}(T^-) \ra \cdots
\eq

\bq\label{khsp}
\cdots \ra \khs^{i-1}(T^+)\ra \khs^{i-1}(TT')\ra \khs^{i-1-2e}(T\cup_i\cap_i T')\ra \khs^{i}(T^+)\ra \cdots 
\eq

\ecor

\vspace{.6cm}

These two exact triangles are the similar to those for Khovanov homology (\cite{knotpolyknothomol}) after the collapse of the bigrading.
The same argument as in \cite{qalternating} can be used to deduce that the two invariants are equivalent  for alternating (and more generally quasi-alternating) links.









\bibliographystyle{plain} 
\bibliography{../biblio}

\def\cprime{$'$}
\begin{thebibliography}{10}

\bibitem{FukayaOh}
Kenji Fukaya and Yong-Geun Oh.
\newblock Zero-loop open strings in the cotangent bundle and {M}orse homotopy.
\newblock {\em Asian J. Math.}, 1(1):96--180, 1997.

\bibitem{categorification}
Mikhail Khovanov.
\newblock A categorification of the {J}ones polynomial.
\newblock {\em Duke Math. J.}, 101(3):359--426, 2000.

\bibitem{functorvalued}
Mikhail Khovanov.
\newblock A functor-valued invariant of tangles.
\newblock {\em Algebr. Geom. Topol.}, 2:665--741 (electronic), 2002.

\bibitem{qalternating}
Ciprian Manolescu and Peter Ozsv{\'a}th.
\newblock On the {K}hovanov and knot {F}loer homologies of quasi-alternating
  links.
\newblock In {\em Proceedings of {G}\"okova {G}eometry-{T}opology {C}onference
  2007}, pages 60--81. G\"okova Geometry/Topology Conference (GGT), G\"okova,
  2008.

\bibitem{PerutzI}
Tim Perutz.
\newblock Lagrangian matching invariants for fibred four-manifolds. {I}.
\newblock {\em Geom. Topol.}, 11:759--828, 2007.

\bibitem{knotpolyknothomol}
Jacob Rasmussen.
\newblock Knot polynomials and knot homologies.
\newblock In {\em Geometry and topology of manifolds}, volume~47 of {\em Fields
  Inst. Commun.}, pages 261--280. Amer. Math. Soc., Providence, RI, 2005.

\bibitem{RR1}
Reza Rezazadegan.
\newblock Seidel-{S}mith cohomology for tangles.
\newblock {\em Selecta Mathematica New Series}, 15:487--518, 2009.

\bibitem{gradedlag}
Paul Seidel.
\newblock Graded {L}agrangian submanifolds.
\newblock {\em Bull. Soc. Math. France}, 128(1):103--149, 2000.

\bibitem{seideltriangle}
Paul Seidel.
\newblock A long exact sequence for symplectic {F}loer cohomology.
\newblock {\em Topology}, 42(5):1003--1063, 2003.

\bibitem{SS}
Paul Seidel and Ivan Smith.
\newblock A link invariant from the symplectic geometry of nilpotent slices.
\newblock {\em Duke Math. J.}, 134(3):453--514, 2006.

\bibitem{SS2}
Paul Seidel and Ivan Smith.
\newblock Localization for involutions in {F}loer cohomology.
\newblock {\em Geom. Funct. Anal.}, 20(6):1464--1501, 2010.

\bibitem{Waldron1}
Jack Waldron.
\newblock An invariant of link cobordisms from symplectic khovanov homology.
\newblock {\em arXiv:0912.5067}.

\bibitem{WWtriangle}
K.~Wehrheim and Woodward C.
\newblock Exact triangle for fibred {D}ehn twists.
\newblock {\em preprint available at math.rutgers.edu/$\sim$ctw}.

\bibitem{WWorient}
Katerin Wehrheim and Christopher Woodward.
\newblock Orientations for pseudoholomorphic quilts.
\newblock {\em To apear}.

\bibitem{WWquilts}
Katerin Wehrheim and Christopher Woodward.
\newblock Pseudoholomorphic quilts.
\newblock {\em ArXiv:0905.1369}.

\bibitem{WW}
Katrin Wehrheim and Christopher~T. Woodward.
\newblock Functoriality for {L}agrangian correspondences.
\newblock {\em arXiv:0708.2851v1}.

\bibitem{Yetter}
David~N. Yetter.
\newblock Markov algebras.
\newblock In {\em Braids (Santa Cruz, CA, 1986)}, volume~78 of {\em Contemp.
  Math.}, pages 705--730. Amer. Math. Soc., Providence, RI, 1988.

\end{thebibliography}



\end{document}